\documentclass[reqno]{amsart}
\usepackage{amssymb}
\usepackage{hyperref}
\usepackage{graphicx,graphpap}
\begin{document}

\title[Existence and decay of solutions]
{Existence and decay of solutions of a nonlinear viscoelastic
problem with a mixed nonhomogeneous condition}

\author[N. T. Long, A. P. N. Dinh, L. X. Truong]
{Nguyen Thanh Long, Alain Pham Ngoc Dinh, Le Xuan Truong}

\address{Nguyen Thanh Long \hfill\break
 Department of Mathematics,
 Hochiminh City National University,
 227 Nguyen Van Cu, Q5, HoChiMinh City, Vietnam}
\email{longnt@hcmc.netnam.vn, longnt2@gmail.com}

\address{Alain Pham Ngoc Dinh \hfill\break
MAPMO, UMR 6628,  Mathematics Dept., University of Orleans, BP 6759,
45067 Orleans Cedex 2, France.} \email{alain.pham@univ-orleans.fr,
alain.pham@math.cnrs.fr}

\address{Le Xuan Truong \hfill\break
Department of Mathematics, Faculty of General Science,
University of Technical Education in HoChiMinh City,
01 Vo Van Ngan Str., Thu Duc Dist., HoChiMinh City, Vietnam}
\email{lxuantruong@gmail.com, truong@math.net}

\keywords{Faedo-Galerkin method; Global existence; contraction \hfill\break\indent mapping theorem; nonlinear wave equation; Viscoelastic; Exponential decay}

\begin{abstract}
We study the initial-boundary value problem for a nonlinear wave equation given by
\begin{gather*}        u_{tt}-u_{xx}+\int_0^{t}k(t-s)u_{xx}(s)ds+|u_{t}|^{q-2}u_{t}=f(x,t,u),\\
u_{x}(0,t)=u(0,t),  u_{x}(1,t)+\eta u(1,t)=g(t),\\
u(x,0)=u_0(x),  u_{t}(x,0)=u_1(x),
\end{gather*}
where $\eta \geq 0$, $q\geq 2$ are given constants and $u_0, u_1, g, k, f$ are given functions.\\
In this paper, we consider two main parts. In Part 1, under a
certain local Lipschitzian condition on $f$ with
$(\widetilde{u}_0,\widetilde{u}_1)\in H^1 \times L^2$; $k,g \in
H^1(0,T)$, $\eta \geq 0$; $q\geq 2$, a global existence and
uniqueness theorem is proved. The proof is based on the paper [10]
associated to a contraction mapping theorem and standard arguments
of density. In Part 2, the asymptotic behavior of the solution $u$
as $t \rightarrow +\infty$ is studied, under more restrictive
conditions, namely $g=0$, $f(x,t,u)=-|u|^{p-2}u +F(x,t)$, $p\geq 2$,
$F \in L^1(\mathbb{R}_+;L^2) \bigcap L^2(\mathbb{R}_+;L^2)$,
$\int_0^{+\infty}e^{\sigma t}\|F(t)\|^2dt<+\infty$, with $\sigma >
0$, and $(\widetilde{u}_0,\widetilde{u}_1)\in H^1\times L^2$, $k\in
H^1(\mathbb{R}_+)$, and some others ($\|\cdot\|$ denotes the
$L^2(0,1)$ norm). It is proved that under these conditions, a unique
solution $u(t)$ exists on $\mathbb{R}_+$ such that $\|u^{/}(t)\|
+\|u_x(t)\|$ decay exponentially to $0$ as $t\rightarrow +\infty$.
Finally, we present some numerical results.
\end{abstract}
\maketitle

\numberwithin{equation}{section}
\newtheorem{theorem}{Theorem}[section]
\newtheorem{lemma}[theorem]{Lemma}
\newtheorem{remark}[theorem]{Remark}


\section{Introduction}
In this paper we will consider the following initial and boundary value problem:
\begin{gather} \label{e1.1}
u_{tt}-u_{xx}+\int_0^{t}k(t-s)u_{xx}(s)ds+|u_{t}|^{q-2}u_{t}=f(x,t,u),0<x<1; 0<t<T,\\
\label{e1.2}u_{x}(0,t)=u(0,t),  u_{x}(1,t)+\eta u(1,t)=g(t),\\
\label{e1.3}u(x,0)=u_0(x),  u_{t}(x,0)=u_1(x),
\end{gather}
where $\eta \geq 0$, $q\geq 2$ are given constants and $u_0, u_1, g, k, f$ are given functions satisfying conditions specified later.

In a recent paper [1], Berrimia and Messaoudi considered the problem
\begin{gather} \label{e1.4}
u_{tt}-\Delta u+\int_0^{t}k(t-s)\Delta u(s)ds=|u|^{p-2}u,  x\in \Omega, t>0,\\
\label{e1.5}u=0, \quad \text{on} \quad \partial\Omega,\\
\label{e1.6}u(x,0)=\widetilde{u}_0(x),u_{t}(x,0)=\widetilde{u}_1(x), x\in \Omega,
\end{gather}
where $p>2$ is a constant, $k$ is a given positive function, and $\Omega$ is a bounded domain of $\mathbb{R}^n$ ($n\geq 1$), with a smooth boundary $\partial \Omega$. This type of problems have been considered by many authors and several results concerning existence, nonexistence, and asymptotic behavior have been established. In this regard, Cavalcanti et al. [3] studied the following equation
\begin{equation}\label{e1.7}
u_{tt}-\Delta u+\int_0^{t}k(t-s)\Delta u(s)ds+|u|^{p-2}u+a(t)u_{t}=0,\quad \text{in} \quad \Omega \times (0,\infty),
\end{equation}
for $a: \Omega \rightarrow \mathbb{R}_+$, a function, which may be null on a part of the domain $\Omega$. Under the conditions that $a(x)\geq a_0>0$ on $\omega \subset \Omega$, with $\omega$ satisfying some geometry restrictions and
\begin{equation}\label{e1.8}
-\zeta_1k(t)=k^{/}(t)=-\zeta_2k(t), t\geq 0,
\end{equation}
the authors established an exponential rate of decay.\\
\indent In [2] Bergounioux, Long and Dinh studied problem \eqref{e1.1}, \eqref{e1.3} with $k=0, q=2,  f(x,t,u)=-Ku+F(x,t)$, and the mixed boundary conditions \eqref{e1.2} standing for
\begin{gather}\label{e1.9}
u_{x}(0,t)=g(t)+hu(0,t)-\int_0^{t}H(t-s)u(0,s)ds,\\
\label{e1.10} u_{x}(1,t)+K_1u(1,t)+\lambda_1u_{t}(1,t)=0,
\end{gather}
where $h\geq 0$, $K, \lambda, K_1, \lambda_1$ are given constants and $g, H$ are given functions.\\
\indent In [7], Long, Dinh and Diem obtained the unique existence, regularity and asymptotic expansion of the problem \eqref{e1.1}, \eqref{e1.3}, \eqref{e1.9} and \eqref{e1.10} in the case of $k=0$, $f(x,t,u)=-K|u|^{p-2}u+F(x,t)$, with $p\geq 2$, $q\geq 2$; $K, \lambda$ are given constants.\\
\indent In [8], Long, Ut and Truc gave the unique existence, stability, regularity in time variable and asymptotic expansion for the solution of problem \eqref{e1.1}- \eqref{e1.3} when $k=0$, $q=2$, $f(x,t,u)=-Ku+F(x,t)$ and $(\widetilde{u}_0,\widetilde{u}_1)\in H^2 \times H^1$. In this case, the problem \eqref{e1.1}- \eqref{e1.3} is the mathematical model describing a shock problem involving a linear viscoelastic bar.\\
\indent In [9], Long and Giai obtained the unique existence and asymptotic expansion for the solution of problem \eqref{e1.1}, \eqref{e1.3} when $k=0$, $q=2$, $f(x,t,u)= -Ku+F(x,t)$ and $(\widetilde{u}_0,\widetilde{u}_1)\in H^1 \times L^2$, and the mixed boundary conditions \eqref{e1.2} standing for
\begin{gather}\label{e1.11}
\begin{aligned}
u_{x}(0,t)=&g(t)+K_1|u(0,t)|^{\alpha-2}u(0,t)+\lambda_1|u_{t}(0,t)|^{\beta-2}u_{t}(0,t)\\
&-\int_0^{t}H(t-s)u(0,s)ds,
\end{aligned}\\ \label{e1.12}
u(1,t)=0,
\end{gather}
where $K$, $\lambda$, $K_1$, $\lambda_1$, $\alpha$, $\beta$ are given constants and $g, H$ are given functions. In this case, the problem \eqref{e1.1}, \eqref{e1.3}, \eqref{e1.11}, \eqref{e1.12} is the mathematical model describing a shock problem involving a nonlinear viscoelastic bar.\\
\indent In [10], Long and Truong obtained the unique existence and asymptotic expansion for the solution of problem \eqref{e1.1} -\eqref{e1.3} when $f(x,t,u)=-K|u|^{p-2}u +F(x,t)$, $(\widetilde{u}_0,\widetilde{u}_1)\in H^2\times H^1$; $F, F_{t}\in L^2(Q_{T})$, $k\in W^{2,1}(0,T)$, $g\in H^2(0,T)$; $K$, $\eta \geq 0$, $\eta_0 > 0$; $p, q\geq 2$.\\
\indent In this paper, we consider two main parts. In Part 1, under
a certain local Lipschitzian condition on $f$ with
$(\widetilde{u}_0,\widetilde{u}_1)\in H^1\times L^2$; $k, g \in
H^1(0,T)$,$\lambda >0$, $\eta_0 >0$; $\eta\geq 0$; $q\geq 2$, a
global existence and uniqueness theorem is proved. The proof is
based on the paper [10] associated to a contraction mapping theorem
and standard arguments of density. In Part 2, the asymptotic
behavior of the solution $u$ as $t\rightarrow \infty$ is studied,
under more restrictive conditions, namely $f(x,t,u)=-|u|^{p-2}u
+F(x,t)$, $p\geq 2$, $F \in L^1(\mathbb{R}_+;L^2) \bigcap
L^2(\mathbb{R}_+;L^2)$, $\int_0^{+\infty}e^{\sigma
t}\|F(t)\|^2dt<+\infty$, with $\sigma>0$, and
$(\widetilde{u}_0,\widetilde{u}_1)\in H^1\times L^2$, $g=0$, $k\in
H^1(\mathbb{R}_+$, and some others ($\|\cdot\|$ denotes the
$L^2(0,1)$ norm). It is proved that under these conditions, a unique
solution $u(t)$ exists on $\mathbb{R}_+$ such that $\|u^{/}(t)\|
+\|u_x(t)\|$ decay exponentially to $0$ as $t\rightarrow +\infty$.
The results obtained here relatively are in part generalizations of
those in [1-3, 6-10]. Finally, we present some numerical results.

\section{Preliminary Results}
Put $\Omega =(0,1)$, $Q_{T}=\Omega \times (0,T)$, $T>0$. We omit
the definitions of usual function spaces: $C^{m}(\overline{\Omega
}) $, $L^{p}(\Omega) $, $W^{m,p}(\Omega) $.  We denote
$W^{m,p}=W^{m,p}(\Omega) $, $ L^{p}=W^{0,p}(\Omega) $,
$H^{m}=W^{m,2}(\Omega) $, $1\leq p\leq \infty $, $m=0,1,\dots $
The norm in $L^2$ is denoted by $\|\cdot \|$. We also denote by
  $\langle \cdot ,\cdot \rangle $
the scalar product in $ L^2$ or pair of dual scalar product of
continuous linear functional with an element of a function space.
We denote by $\|\cdot \|_{X}$ the norm of a Banach space $X$ and by $X'$ the dual space of $X$. We denote by $L^{p}(0,T;X)$, $1\leq p\leq \infty $ for the Banach space of the
real functions $u:(0,T)\rightarrow X$ measurable, such that
$$
\|u\|_{L^{p}(0,T;X)}=\Big(\int_0^{T}\|
u(t)\|_{X}^{p}dt\Big) ^{1/p}<\infty \quad \text{for } 1\leq p<\infty,
$$
and
$$
\|u\|_{L^{\infty }(0,T;X)}= \mathop{\rm ess\,sup}_{0<t<T}
\|u(t)\|_{X}\quad\text{for }p=\infty .
$$
Let $u(t)$, $u'(t)=u_t(t)$, $u''(t)=u_{tt}(t)$, $u_{x}(t)$, and
$u_{xx}(t)$ denote $u(x,t)$, $\frac{\partial u}{\partial t}(x,t)$,
$\frac{\partial ^2u}{\partial t^2}(x,t)$,
$\frac{\partial u}{\partial x}(x,t)$, and
$\frac{\partial ^2u}{\partial x^2}(x,t)$, respectively.
\vspace{0.5cm}\\
Without loss of generality, we can suppose that $\eta_0=\lambda= 1$. For every $\eta \geq 0$, we put

\begin{gather}\label{e2.1}
a_{\eta}(u,v)=\int_0^1u_{x}(x)v_{x}(x)dx+u(0)v(0)+\eta u(1)v(1),  \forall u,v\in H^1,\\
\label{e2.2} \|v\|_{\eta}=(a_{\eta}(v,v))^{1/2}.
\end{gather}
On $H^1$ we shall use the following equivalent norm
\begin{equation}\label{e2.3}
\|v\|_1=\left(v^2(0)+\int_0^1|v_{x}(x)|^2dx\right)^{1/2}
\end{equation}
Then we have the following lemmas.

\begin{lemma} \label{lem2.1}
The imbedding $V\hookrightarrow C^{0}([0,1])$ is compact and
\begin{equation} \label{e2.4}
\|v\|_{C^{0}([0,1])}\leq \|v\|_{V}, \text{ for all } \ v\in V.
\end{equation}
\end{lemma}
\begin{lemma} \label{lem2.2}
Let $\eta\geq 0$. Then, the symmetric bilinear form $a_{\eta}(\cdot,\cdot)$ defined by \eqref{e2.1} is continuous on $H^1 \times H^1$ and coercive on $H^1$, i.e.,\\
(i) $|a_{\eta}(u,v)|=C_{\eta}\|u\|_1\|v\|_1,$ \quad     \text{for all} \quad  $u, v \in H^1$, \\
(ii) $a_{\eta}(v,v)=\|v\|_1^2$, \quad \text{for all} \quad $v \in H^1$,\\
where $C_{\eta}=1+2\eta.$
\end{lemma}
\noindent
The proofs of these lemmas are straightforward, and we omit the details.\\
We also note that on $H^1$, $\|v\|_1, \|v\|_{H^1} = \left(\|v\|^2+\|v^{/}\|^2\right)^{1/2}$, $\|v\|_{\eta}=\sqrt{(a_{\eta}(v,v))}$ are three equivalent norms.
\begin{gather}\label{e2.5}
\|v\|_1^2 \leq \|v\|_{\eta}^2\leq C_{\eta}\|v\|_1^2, \quad \text{for all} \quad v\in H^1,\\
\label{e2.6} \frac{1}{3}\|v\|_{H^1}^2\leq \|v\|_1^2\leq 3\|v\|_{H^1}^2, \quad \text{for all} \quad v\in H^1,
\end{gather}

\section{ The Existence and uniqueness theorem of the solution}
In this section we study the global existence of solutions for problem \eqref{e1.1}-\eqref{e1.3}. For this purpose, we consider, first, a related nonlinear problem. Then, we use the well-known Banach's fixed point theorem to prove the existence of solutions to the nonlinear problem \eqref{e1.1}-\eqref{e1.3}.\\
\indent We make the following assumptions:
\begin{itemize}
\item[(H1)] $\eta \geq 0, q\geq 2$,
\item[(H2)] $k, g\in H^1(0,T)$,
\item[(H3)] $\widetilde{u}_0\in H^1$ and $\widetilde{u}_1\in L^2,$
\item[(H4)] $f\in C^0(\overline{\Omega}\times \mathbb{R}_+\times \mathbb{R})$ satisfies the conditions $D_2f, D_3f \in C^0(\overline{\Omega}\times \mathbb{R}_+\times \mathbb{R}).$
\end{itemize}
\vspace{0.15cm}
For each $T>0$, we put

\begin{equation}\label{e3.1}
W(T)={v\in L^{\infty}(0,T;H^1):v_{t}\in L^{\infty}(0,T;L^2)\bigcap L^{q}(Q_{T})}.
\end{equation}
Then $W(T)$ is a Banach space with respect to the norm (see[5]):
\begin{equation}\label{e3.2}
\|v\|_{W(T)}=\|v\|_{L^{\infty}(0,T;H^1)}+\|v_{t}\|_{L^{\infty}(0,T;L^2)}+\|v_{t}\|_{L^{q}(Q_{T})},  v\in W(T).
\end{equation}
For each $v\in W(T)$, we associate with the problem \eqref{e1.1}-\eqref{e1.3} the following variational problem.\\
Find $u\in W(T)$ which satisfies the variational problem
\begin{equation}\label{e3.3}
\begin{aligned}
&<u^{//}(t),w>+a_{\eta}(u(t),w)-\int_0^{t}k(t-s)a_{\eta}(u(s),w)ds +<\psi_{q}(u^{/}(t)),w>\\
&=g_1(t)w(1)+<f(\cdot,t,v(\cdot,t)),w> \quad \text{for all} \quad w\in H^1,
\end{aligned}
\end{equation}
\begin{equation}\label{e3.4}
u(0)=\widetilde{u}_0, u_{t}(0)=\widetilde{u}_1,
\end{equation}
where
\begin{equation}\label{e3.5}
\psi_{q}(z)=|z|^{q-2}z, g_1(t)=g(t)-\int_0^{t}k(t-s)g(s)ds.
\end{equation}
Then, we have the following theorem
\begin{theorem} \label{thm3.1}
Let (H1)-(H4) hold. Then, for every $T>0$ and $v\in W(T)$, problem \eqref{e3.3}- \eqref{e3.5} has a unique solution $u\in W(T)$ and such that
\begin{equation}\label{e3.6}
u^{//}, u_{xx}\in L^{q^{/}}(0,T;(H^1)^{/}), \quad \text{where} \quad
q^{/}=q/(q-1).
\end{equation}
Furthermore, we have
\begin{equation}\label{e3.7}
\|u^{/}(t)\|^2+\|u(t)\|_{\eta}^2+2 \int_0^{t} \|u^{/}(s)\|_{L^{q}}^{q}ds \leq C_{1T}\exp(TC_{2T}), \forall t\in [0,T],
\end{equation}
where
\begin{equation}\label{e3.8}
\begin{gathered}
C_{1T}=C_{1T}(v,\widetilde{u}_0,\widetilde{u}_1,k,g)= 2\left[\|\widetilde{u}_1\|^2 + \|\widetilde{u}_0\|_{\eta}^2 +2|g_1(0)\widetilde{u}_0(1)|
\begin{array}{c}
  \\
  \\
\end{array} \right. \\
\left. + 6\|g_1\|_{L^{\infty}(0,T)}^2 +2\|g_1^{/}\|_{L^2(0,T)}^2 +\int_0^{T}\|f(\cdot,s,v(s))\|^2ds\right],
\end{gathered}
\end{equation}
\begin{equation}\label{e3.9}
C_{2T}=C_{2T}(k)=2\left[3+2|k(0)|+6\|k\|_{L^2(0,T)}^2+T\|k^{/}\|_{L^2(0,T)}^2\right],
\end{equation}
and
\begin{equation}\label{e3.10}
g_1(t)=g(t)-\int_0^{t}k(t-s)g(s)ds.
\end{equation}
\end{theorem}
\begin{proof}[Proof of theorem 3.1]
The proof consists of steps two steps

\noindent {\textbf{a. The existence of solution.}}
We approximate $\widetilde{u}_0$, $\widetilde{u}_1$, $k$, $g$ by sequences $\{u_{0m}\}\subset C_0^{\infty} \left(\overline{\Omega}\right)$, ${u_{1m}}\subset C_0^{\infty}(\Omega)$, ${k_{m}}, {g_{m}}\subset C_0^{\infty}([0,T])$, respectively, such that
\begin{equation}\label{e3.11}
\begin{gathered}
u_{0m}\rightarrow \widetilde{u}_0 \quad \text{strongly in} \quad H^1, \\
u_{1m}\rightarrow \widetilde{u}_1 \quad \text{strongly in} \quad L^2,\\
k_{m}\rightarrow k \quad \text{strongly in} \quad H^1(0,T),  \\
g_{m}\rightarrow g \quad \text{strongly in} \quad H^1(0,T).
\end{gathered}
\end{equation}
Then we consider the following variational problem: Find $u_{m}\in W(T)$ which satisfies the variational problem
\begin{equation}\label{e3.12}
\begin{aligned}
&<u_{m}^{//}(t),w>+a_{\eta}(u_{m}(t),w)-\int_0^{t}k_{m}(t-s)a_{\eta}(u_{m}(s),w)ds\\
&+<\psi_{q}(u_{m}^{/}(t)),w>=g_{1m}(t)w(1)+<f(\cdot,t,v(\cdot,t)),w>, \forall w\in H^1,
\end{aligned}
\end{equation}
\vspace{0.2cm}
\begin{equation}\label{e3.13}
u(0)=u_{0m},  u^{/}(0)=u_{1m},
\end{equation}
and
\begin{equation}\label{e3.14}
u_{m}\in L^{\infty}(0,T;H^2), u_{m}^{/}\in L^{\infty}(0,T;H^1), u_{m}^{//}\in L^{\infty}(0,T;L^2),
\end{equation}
where
\begin{equation}\label{e3.15}
g_{1m}(t)=g_{m}(t)-\int_0^{t}k_{m}(t-s)g_{m}(s)ds.
\end{equation}
The existence of a sequence of solutions ${u_{m}}$ satisfying \eqref{e3.12}-\eqref{e3.15} is a direct result of the theorem 2.1 in [10]. We shall prove that ${u_{m}}$ is a Cauchy sequence in $W(T)$.\\
{\textit{(i) A priori estimates.}}\\
We take $w=u_{m}^{/}(t)$ in \eqref{e3.12}, afterwards integrating with respect to the time variable from $0$ to $t$, we get after some rearrangements

\begin{equation}\label{e3.16}
\begin{aligned}
\sigma_{m}(t)=&\sigma_{m}(0)-2g_{1m}(0)u_{0m}(1)+2g_{1m}(t)u_{m}(1,t)\\
&-2\int_0^{t}g_{1m}^{/}(r)u_{m}(1,r)dr -2k_{m}(0)\int_0^{t}\|u_{m}(r)\|_{\eta}^2dr \\
&+2\int_0^{t}k_{m}(t-s)a_{\eta}(u_{m}(s),u_{m}(t))ds\\ &-2\int_0^{t}dr\int_0^{r}k_{m}^{/}(r-s)a_{\eta}(u_{m}(s),u_{m}(r))ds\\
&+2\int_0^{t}<f(\cdot,s,v(\cdot,s)),u_{m}^{/}(s)>ds,
\end{aligned}
\end{equation}
where
\begin{equation}\label{e3.17}
\sigma_{m}(t)=\|u_{m}^{/}(t)\|^2+\|u_{m}(t)\|_{\eta}^2+2 \int_0^{t}||u_{m}^{/}(s)||_{L^{q}}^{q}ds.
\end{equation}
Proving in the same manner as in [10], we have the following results:
\begin{equation}\label{e3.18}
\sigma_{m}(t)=C_{1T}(m)+C_{2T}(m)\int_0^{t}\sigma_{m}(s)ds,
\forall t\in [0,T],
\end{equation}
where
\begin{gather}\label{e3.19}
\begin{aligned}
C_{1T}(m)=&2\left[ \|u_{1m}\|^2+\|u_{0m}\|_{\eta}^2+2|g_{1m}(0)u_{0m}(1)|+6\|g_{1m}\|_{L^{\infty}(0,T)}^2\right.\\
&\left.+2\|g_{1m}^{/}\|_{L^2(0,T)}^2+\int_0^{T}\|f(\cdot,s,v(s))\|^2ds\right],
\end{aligned}\\
\label{e3.20}
C_{2T}(m)=2\left[3+2|k_{m}(0)|+6\|k_{m}\|_{L^2(0,T)}^2+T\|k_{m}^{/}\|_{L^2(0,T)}^2\right].
\end{gather}
From the assumptions (H1)-(H4), afterwards using Gronwall's lemma, we deduce from \eqref{e3.11}, that
\begin{equation}\label{e3.21}
\sigma_{m}(t)\leq \widetilde{C}_{T}, \quad \text{for all} \quad m \quad \text{and} \quad t\in [0,T],
\end{equation}
where $C_{T}$ is a constant independent of $m$.

On the other hand, we deduce from \eqref{e3.12}, \eqref{e3.21}, that, for  all $w\in H^1$, we have
\begin{equation}\label{e3.22}
\begin{aligned}
|<u_{m}^{//}(t),w>|\leq &\|u_{m}(t)\|_{\eta}\|w\|_{\eta} +\int_0^{t}|k_{m}(t-s)|\|u_{m}(s)\|_{\eta}\|w\|_{\eta}ds\\
&+\|\psi_{q}(u_{m}^{/})\|_{L^{q^{/}}(\Omega)}\|w\|_{L^{q}(\Omega)}+|g_{1m}(t)|\|w\|_{\eta}\\
&+\|f(\cdot,t,v(\cdot,t))\|\|w\|\\
& \leq C_{T}\sqrt{(3C_{\eta})}\left[1+\|\psi_{q}(u_{m}^{/})\|_{L^{q^{/}}(\Omega)}\right]\|w\|_{H^1}.
\end{aligned}
\end{equation}
This implies that
\begin{equation}\label{e3.23}
\begin{aligned}
\|u_{m}^{//}(t)\|_{(H^1)^{/}}&= \sup_{0\neq w \in H^1}\frac{\left|<u_{m}^{//}(t),w>\right|}{\|w\|_{H^1}}\\
&\leq C_{T}\sqrt{3C_{\eta}}\left[1+\|\psi_{q}(u_{m}^{/})\|_{L^{q^{/}}(\Omega)}\right].
\end{aligned}
\end{equation}
Hence
\begin{equation}\label{e3.24}
\begin{aligned}
\|u_{m}^{//}\|_{L^{q^{/}}(0,T;(H^1)^{/})}^{q^{/}}&=\int_0^{T}\|u_{m}^{//}(t)\|_{(H^1)^{/}}^{q^{/}} dt \\
&\leq\left(C_{T}\sqrt{3C_{\eta}}\right)^{q^{/}}2^{q^{/}-1}\int_0^{T}\left[1+\|u_{m}^{/}(t)\|_{L^{q}(\Omega)}^{q}\right]dt\\
&\leq C_{T},
\end{aligned}
\end{equation}
where $C_{T}$ always indicating a constant depending on $T$.\\
{\textit{(ii) The convergence of sequence $\{u_{m}\}$}}\\
We shall prove that ${u_{m}}$ is a Cauchy sequence in $W(T)$. Let $\widehat{u}=u_{m}-u_{\mu}$. Then $\widehat{u}$  satisfies the variational problem

\begin{equation}\label{e3.25}
\begin{aligned}
&<u^{//}(t),w>+a_{\eta}(\widehat{u}(t),w)-\int_0^{t}k_{m}(t-s)a_{\eta}(\widehat{u}(s),w)ds\\
&-\int_0^{t}\widehat{k}(t-s)a_{\eta}(u_{\mu}(s),w)ds
+<\psi_{q}(u_{m}^{/}(t))-\psi_{q}(u_{\mu}^{/}(t)),w>\\
&= g_1(t)w(1) \quad \text{for all} \quad w\in H^1,
\end{aligned}
\end{equation}
\vspace{0.2cm}
\begin{equation}\label{e3.26}
\widehat{u}(0)=\widehat{u}_0,  \widehat{u}^{/}(0)=\widehat{u}1,
\end{equation}
where
\begin{equation}\label{e3.27}
\begin{gathered}
\widehat{u}_0= u_{0m}-u_{0\mu},  \widehat{u}_1=u_{1m}-u_{1\mu}, \\
\widehat{k}=k_{m}-k_{\mu}, \widehat{g}=g_{m}-g_{\mu}, \widehat{g}_1=g_{1m}-g_{1\mu},\\
\widehat{g}_1(t)= \widehat{g}(t)- \int_0^{t}k_{m}(t-s)\widehat{g}(s)ds-\int_0^{t}\widehat{k}(t-s)g_{\mu}(s)ds.
\end{gathered}
\end{equation}
We take $w=u^{/}(t)$ in \eqref{e3.25}, after integrating with respect to the time variable from $0$ to $t$, we get after some rearrangements

\begin{equation}\label{e3.28}
\begin{aligned}
Z(t)=&Z(0)-2\widehat{g}_1(0)\widehat{u}_0(1)+2\widehat{g}_1(t)\widehat{u}(1,t)-2\int_0^{t}\widehat{g}1^{/}(r)\widehat{u}(1,r)dr \\
&-2k_{m}(0)\int_0^{t}\|\widehat{u}(r)\|_{\eta}^2dr+2\int_0^{t}k_{m}(t-s)a_{\eta}(\widehat{u}(s),\widehat{u}(t))ds \\
&-2\int_0^{t}dr\int_0^{r}k_{m}^{/}(r-s)a_{\eta}(\widehat{u}(s),\widehat{u}(r))ds \\
&-2\widehat{k}(0)\int_0^{t}a_{\eta}(u_{\mu}(s),\widehat{u}(s))ds +2\int_0^{t}\widehat{k}(t-s)a_{\eta}(u_{\mu}(s),\widehat{u}(t))ds\\
&-2\int_0^{t}dr\int_0^{r}\widehat{k}^{/}(r-s)a_{\eta}(u_{\mu}(s),\widehat{u}(r))ds,
\end{aligned}
\end{equation}
where
\begin{equation}\label{e3.29}
\begin{aligned}
Z(t)=&\|\widehat{u}^{/}(t)\|^2+\|\widehat{u}(t)\|_{\eta}^2\\
&+2\int_0^{t}<\psi_{q}(u_{m}^{/}(s))-\psi_{q}(u_{\mu}^{/}(s)),u_{m}^{/}(s)-u_{\mu}^{/}(s)>ds.
\end{aligned}
\end{equation}
Using the following inequality
\begin{equation}\label{e3.30}
\forall q\geq 2, \exists C_{q}>0: (|x|^{q-2}x-|y|^{q-2}y)(x-y)\geq C_{q}|x-y|^{q},  \forall x, y \in \mathbb{R},
\end{equation}
it follows from \eqref{e3.29} that
\begin{equation}\label{e3.31}
Z(t)\geq \|\widehat{u}^{/}(t)\|^2+\|\widehat{u}(t)\|_{\eta}^2+2C_{q}\int_0^{t}\|\widehat{u}^{/}(s)\|_{L^{q}}^{q}ds.
\end{equation}
Using the inequality
\begin{equation}\label{e3.32}
2ab \leq \epsilon a^2+\frac{1}{\epsilon}b^2, \forall a,b\in \mathbb{R}, \forall \epsilon >0,
\end{equation}
and the following inequalities
\begin{equation}\label{e3.33}
|a_{\eta}(u,v)|\leq \|u\|_{\eta}\|v\|_{\eta}, \forall u,v \in H^1,
\end{equation}
\begin{equation}\label{e3.34}
|\widehat{u}(1,t)|\leq \|\widehat{u}(t)\|_{C^0(\Omega)}\leq \sqrt{2}\|\widehat{u}(t)\|_1\leq \sqrt{2}\|\widehat{u}(t)\|_{\eta}\leq \sqrt{2Z(t)},
\end{equation}
 we shall estimate respectively the following terms on the right-hand side of \eqref{e3.28} as follows
\begin{equation}\label{e3.35}
\begin{aligned}
Z(0)-2\widehat{g}_1(0)\widehat{u}_0(1)\leq & \|u_{1m}-u_{1\mu}\|^2+\|u_{0m}-u_{0\mu}\|_{\eta}^2  \\
&+2|g_{1m}(0)-g_{1\mu}(0)| |u_{0m}(1)-u_{0\mu}(1)|,
\end{aligned}
\end{equation}
\begin{equation}\label{e3.36}
2\widehat{g}_1(t)\widehat{u}(1,t)\leq 8\|\widehat{g}_1\|_{L^{\infty}(0,T)}^2 +\frac{1}{4}Z(t), \quad \text{with} \quad \epsilon =\frac{1}{8},
\end{equation}
\begin{equation}\label{e3.37}
-2\int_0^{t}\widehat{g}_1^{/}(r)\widehat{u}(1,r)dr\leq 2\|\widehat{g}_1^{/}\|_{L^2(0,T)}^2+\int_0^{t}Z(r)dr,
\end{equation}
\begin{equation}\label{e3.38}
2\int_0^{t}k_{m}(t-s)a_{\eta}(\widehat{u}(s),\widehat{u}(t))ds\leq \frac{1}{8}Z(t)+8\|k_{m}\|_{L^2(0,T)}^2\int_0^{t}Z(s) ds,
\end{equation}
\begin{equation}\label{e3.39}
-2k_{m}(0)\int_0^{t}\|\widehat{u}(r)\|_{\eta}^2dr\leq  2|k_{m}(0)|\int_0^{t}Z(r)dr,
\end{equation}
\begin{equation}\label{e3.40}
-2\int_0^{t}dr\int_0^{r}k_{m}^{/}(r-s)a_{\eta}(\widehat{u}(s),\widehat{u}(r))ds\leq
\left(1+T\|k_{m}^{/}\|_{L^2(0,T)}^2\right)\int_0^{t}Z(s) ds,
\end{equation}
\begin{equation}\label{e3.41}
2\int_0^{t}\widehat{k}(t-s)a_{\eta}(u_{\mu}(s),\widehat{u}(t))ds\leq
\frac{1}{8}Z(t)+8\widetilde{C}_{T}\|\widehat{k}\|_{L^1(0,T)}^2,
\end{equation}
\begin{equation}\label{e3.42}
-2\widehat{k}(0)\int_0^{t}a_{\eta}(u_{\mu}(s),\widehat{u}(s))ds\leq
T\widetilde{C}_{T}|\widehat{k}(0)|^2+\int_0^{t}Z(s) ds,
\end{equation}
\begin{equation}\label{e3.43}
-2\int_0^{t}dr\int_0^{r}\widehat{k}^{/}(r-s)a_{\eta}(u_{\mu}(s),\widehat{u}(r))ds\leq
T^2\widetilde{C}_{T}\|\widehat{k}^{/}\|_{L^2(0,T)}^2+\int_0^{t}Z(s)ds.
\end{equation}
Combining \eqref{e3.28}, \eqref{e3.29}, \eqref{e3.31} and
\eqref{e3.35}-\eqref{e3.43}, we obtain
\begin{equation}\label{e3.44}
Z(t)\leq \rho_{1T}(m,\mu)+\rho_{2T}(m)\int_0^{t}Z(s)ds, \forall t
\in [0,T],
\end{equation}
where
\begin{equation}\label{e3.45}
\begin{gathered}
\rho_{1T}(m,\mu)=2\left[\|\widehat{u}_1\|^2+\|\widehat{u}_0\|_{\eta}^2+2|\widehat{g}_1(0)\widehat{u}_0(1)|
+8\|\widehat{g}_1\|_{L^{\infty}(0,T)}^2\right.\\
\left.
+2\|\widehat{g}_1^{/}\|_{L^2(0,T)}^2+8\widetilde{C}_{T}\|\widehat{k}\|_{L^1(0,T)}^2+T\widetilde{C}_{T}
|\widehat{k}(0)|^2+T^2\widetilde{C}_{T}\|\widehat{k}^{/}\|_{L^2(0,T)}^2\right],\\
\rho_{2T}(m)=2\left[4+2|k_{m}(0)|+8\|k_{m}\|_{L^2(0,T)}^2+T\|k_{m}^{/}\|_{L^2(0,T)}^2\right].
\end{gathered}
\end{equation}
By Gronwall's lemma, we deduce from \eqref{e3.31}, \eqref{e3.44},
\eqref{e3.45}, that
\begin{equation}\label{e3.46}
\begin{aligned}
\|\widehat{u}^{/}(t)\|^2&+\|\widehat{u}(t)\|_{\eta}^2+2C_{q}\int_0^{t}\|u^{/}(s)\|_{L^{q}}^{q}ds\leq
Z(t)\\
&\leq \rho_{1T}(m,\mu)exp(T\rho_{2T}(m)), \quad \text{for all}
\quad t\in [0,T].
\end{aligned}
\end{equation}
By \eqref{e3.11}, \eqref{e3.27} and \eqref{e3.45}, we obtain
$\rho_{1T}(m,\mu)exp(T\rho_{2T}(m))\rightarrow 0$ as $m, \mu
\rightarrow +\infty$. Hence, it follows from \eqref{e3.46} that
$\{u_{m}\}$ is a Cauchy sequence in $W(T)$. Therefore there exists
$u\in W(T)$ such that
\begin{equation}\label{e3.47}
u_{m}\rightarrow u  \quad \text{strongly  in} \quad  W(T).
\end{equation}
On the other hand, by \eqref{e3.47} and the continuity of
$\psi_{q}$, we obtain
\begin{equation}\label{e3.48}
\psi_{q}(u_{m}^{/})\rightarrow \psi_{q}(u^{/}) \quad \text{a.e.}
\quad (x,t)\in Q_{T}.
\end{equation}
By means of \eqref{e3.21}, it follows that
\begin{equation}\label{e3.49}
\|\psi_{q}(u_{m}^{/})\|_{L^{q^{/}}(Q_{T})}^{q^{/}}=\|u_{m}^{/}\|_{L^{q}(Q_{T})}^{q}\leq
\frac{1}{2}\widetilde{C}_{T},
\end{equation}
for all $m$. By Lions's lemma [5, Lemma 1.3, p. 12], it follows
from \eqref{e3.48} and \eqref{e3.49} that
\begin{equation}\label{e3.50}
\psi_{q}(u_{m}^{/})\rightarrow \psi_{q}(u^{/}) \quad \text{in}
\quad L^{q^{/}}(Q_{T}) \quad \text{weakly}.
\end{equation}
Noticing \eqref{e3.11}$_3$ and \eqref{e3.47} we have
\begin{equation}\label{e3.51}
\begin{aligned}
&\left|\int_0^{T}dt\int_0^{t}k_{m}(t-s)a_{\eta}(u_{m}(s),w(t))ds\right.\\
&\left.-\int_0^{T}dt\int_0^{t}k(t-s)a_{\eta}(u(s),w(t))ds\right|\\
&\leq
\left|\int_0^{T}dt\int_0^{t}k_{m}(t-s)a_{\eta}(u_{m}(s)-u(s),w(t))ds\right|\\
&+\left|\int_0^{T}dt\int_0^{t}[k_{m}(t-s)-k(t-s)]a_{\eta}(u(s),w(t))ds\right|\\
&\leq
3C_{\eta}\|w\|_{L^1(0,T;H^1)}\left[\|k_{m}\|_{L^1(0,T)}\|u_{m}-u\|_{L^{\infty}(0,T;H^1)}\right.\\
&\left.+\|k_{m}-k\|_{L^1(0,T)}\|u\|_{L^{\infty}(0,T;H^1)}\right]
\rightarrow 0
\end{aligned}
\end{equation}
for all $w\in L^1(0,T;H^1)$.

On the other hand, by \eqref{e3.11}$_{3,4}$ and \eqref{e3.15}, we
also obtain
\begin{equation}\label{e3.52}
g_{1m}\rightarrow g_1 \quad \text{strongly  in} \quad H^1(0,T).
\end{equation}
From \eqref{e3.24}, we deduce the existence of a subsequence of
$\{u_{m}\}$, still denoted by $\{u_{m}\}$, such that
\begin{equation}\label{e3.53}
u_{m}^{//}\rightarrow u^{//} \quad \text{in} \quad
L^{q^{/}}(0,T;(H^1)^{/}) \quad \text{weak}.
\end{equation}
Passing to the limit in \eqref{e3.12}, \eqref{e3.13} by
\eqref{e3.47} and \eqref{e3.50}-\eqref{e3.53} we have $u$
satisfying the equation
\begin{equation}\label{e3.54}
\begin{aligned}
&\left<u^{//}(t),w\right>+a_{\eta}(u(t),w)-\int_0^{t}k(t-s)a_{\eta}(u(s),w)ds+\left<\psi_{q}(u^{/}(t)),w\right>\\
&=g_1(t)w(1)+<f(\cdot,t,v(\cdot,t)),w>, \forall w\in H^1, \quad
\text{in} \quad L^{q^{/}}(0,T) \quad \text{weak},
\end{aligned}
\end{equation}
and
\begin{equation}\label{e3.55}
u(0)=\widetilde{u}_0,  u^{/}(0)=\widetilde{u}_1.
\end{equation}
On the other hand, we deduce from \eqref{e3.54}, that
\begin{equation}\label{e3.56}
u_{xx}(t)-\int_0^{t}k(t-s)u_{xx}(s)ds =\phi(t),
\end{equation}
where
\begin{equation}\label{e3.57}
\phi(t)=u^{//}(t)+|u^{/}|^{q-2}u^{/}-f(\cdot,t,v(\cdot,t))\in
L^{q^{/}}(0,T;(H^1)^{/}).
\end{equation}
Hence, it follows from \eqref{e3.56} and \eqref{e3.57}, that
\begin{equation}\label{e3.58}
\begin{aligned}
&\|u_{xx}(t)\|_{(H^1)^{/}}^{q^{/}}\leq
\left(\|\phi(t)\|_{(H^1)^{/}}+\int_0^{t}|k(t-s)|\|u_{xx}(s)\|_{(H^1)^{/}}
ds\right)^{q^{/}}\\
&\leq
2^{q^{/}-1}\left[\|\phi(t)\|_{(H^1)^{/}}^{q^{/}}+\left(\int_0^{t}|k(t-s)|\|u_{xx}(s)\|_{(H^1)^{/}}ds\right)^{q^{/}}\right]\\
&\leq
2^{q^{/}-1}\left[\|\phi(t)\|_{(H^1)^{/}}^{q^{/}}+\|k\|_{L^{q}(0,T)}^{q^{/}}\left(\int_0^{t}\|u_{xx}(s)\|_{(H^1)^{/}}^{q^{/}}ds\right)\right].\\
\end{aligned}
\end{equation}
This implies that
\begin{equation}\label{e3.59}
\begin{aligned}
\int_0^{r}\|u_{xx}(t)\|_{(H^1)^{/}}^{q^{/}}dt\leq &
2^{q^{/}-1}\|\phi\|_{L^{q^{/}}(0,T;(H^1)^{/})}^{q^{/}}\\
&+2^{q^{/}-1}\|k\|_{L^{q}(0,T)}^{q^{/}}\int_0^{r}dt\int_0^{t}\|u_{xx}(s)\|_{(H^1)^{/}}^{q^{/}}
ds
\end{aligned}
\end{equation}
Using Gronwall's lemma, we obtain
\begin{equation}\label{e3.60}
\begin{aligned}
\int_0^{r}\|u_{xx}(t)\|_{(H^1)^{/}}^{q^{/}}dt \leq
2^{q^{/}-1}\|\phi\|_{L^{q^{/}}(0,T;(H^1)^{/})}^{q^{/}}exp\left(2^{q^{/}-1}\|k\|_{L^{q}(0,T)}^{q^{/}}r\right)\leq
C_{T},
\end{aligned}
\end{equation}
where $C_{T}$ always indicating a constant depending on $T$\\
Thus
\begin{equation}\label{e3.61}
u_{xx}\in L^{q^{/}}(0,T;(H^1)^{/}).
\end{equation}
On ther other hand, the estimate \eqref{e3.7} hold by means of
\eqref{e3.11}, \eqref{e3.18}, \eqref{e3.19}, \eqref{e3.20},
\eqref{e3.47}. The existence of the theorem is proved
completely.
\vspace{0.2cm}\\
{\textbf{b. Uniqueness of the solution.}} First, we shall now
require the following lemma.
\begin{lemma}\label{lem3.2}
Let $u$ be the weak solution of  the following problem
\begin{equation}\label{e3.62}
\begin{aligned}
&u^{//}-u_{xx}+\int_0^{t}k(t-s)u_{xx}(s) ds=\Phi, 0<x<1, 0<t<T,\\
&u_{x}(0,t)=u(0,t),   u_{x}(1,t)+\eta u(1,t)=0,\\
&u(x,0)=\widetilde{u}_0(x), u^{/}(x,0)=\widetilde{u}_1(x),\\
&u\in L^{\infty}(0,T;H^1),  u^{/}?L^{\infty}(0,T;L^2),\\
&k\in H^1(0,T),  \Phi\in L^2(Q_{T}).
\end{aligned}
\end{equation}
Then we have
\begin{equation}\label{e3.63}
\begin{aligned}
&\frac{1}{2}\|u^{/}(t)\|^2+\frac{1}{2}\|u(t)\|_{\eta}^2=\frac{1}{2}\|u_1\|^2+\frac{1}{2}\|u_0\|_{\eta}^2
-k(0)\int_0^{t}\|u(r)\|_{\eta}^2dr\\
&+\int_0^{t}k(t-s)a(u(s),u(t))ds-\int_0^{t}dr\int_0^{r}k^{/}(r-s)a(u(s),u(r))ds\\
&+\int_0^{t}<\Phi(s),u^{/}(s)>ds, \quad \text{a.e.} \quad t\in
[0,T].
\end{aligned}
\end{equation}
Furthermore, if  $u_0=u_1=0$ there is equality in \eqref{e3.63}.
\end{lemma}
\noindent The idea of the proof is the same as in [4, Lemma 2.1,
p. 79].
\vspace{0.2cm} \\
We now return to the proof of the
uniqueness of a solution of the problem \eqref{e3.3}-\eqref{e3.5}.
Let $u_1$, $u_2$ be two weak solutions of problem
\eqref{e3.3}-\eqref{e3.5}, such that
\begin{equation}\label{e3.64}
u_{i}\in W(T),  u_{i}^{//},  u_{ixx}\in L^{q^{/}}(0,T;(H^1)^{/}),
i=1, 2.
\end{equation}
Then $u=u_1-u_2$ is the weak solution of the following problem
\begin{equation}\label{e3.65}
\begin{aligned}
&u^{//}-u_{xx}+\int_0^{t}k(t-s)u_{xx}(s)ds+\psi_{q}(u_1^{/})-\psi_{q}(u_2^{/})=0,\\
&u_{x}(0,t)-u(0,t)=u_{x}(1,t)+\eta u(1,t)=0,\\
&u(0)=u^{/}(0)=0,\\
&u\in W(T),  u^{//},  u_{xx}\in L^{q^{/}}(0,T;(H^1)^{/}).
\end{aligned}
\end{equation}
By using Lemma 3.2 with $u_0=u_1=0$,
$\Phi=-\psi_{q}(u_1^{/})+\psi_{q}(u_2^{/})$, we have
\begin{equation}\label{e3.66}
\begin{aligned}
\sigma(t)=&2\int_0^{t}k(t-s)a(u(s),u(t))ds-2k(0)\int_0^{t}\|u(r)\|_{\eta}^2dr\\
&-2 \int_0^{t}dr\int_0^{r}k^{/}(r-s)a(u(s),u(r))ds,
\end{aligned}
\end{equation}
where
\begin{equation}\label{e3.67}
\sigma(t)=\|u^{/}(t)\|^2+\|u(t)\|_{\eta}^2+2\int_0^{t}\left<\psi_{q}(u_1^{/}(s))-\psi_{q}(u_2^{/}(s)),u^{/}(s)\right>ds.
\end{equation}
By using the same computations as in the above part we obtain from
\eqref{e3.66} that
\begin{equation}\label{e3.68}
\sigma(t)=2\left(1+2\|k\|_{L^2(0,T)}^2+2|k(0)|+\|k^{/}\|_{L^1(0,T)}^2\right)\int_0^{t}\sigma(r)dr.
\end{equation}
By Gronwall's lemma, we deduce that $\sigma(t)=0$ and Theorem 3.1
is completely proved.
\end{proof}
\begin{theorem}\label{thm3.3}
Let $T>0$ and $(H1)-(H4)$ hold. Then there exists $T_1\in (0,T)$
such that problem \eqref{e1.1}- \eqref{e1.3} has a unique weak
solution $u\in W(T_1)$ and such that
\begin{equation}\label{e3.69}
u^{//},  u_{xx}\in L^{q^{/}}(0,T_1;(H^1)^{/}).
\end{equation}
\end{theorem}
\begin{proof}
For each $T_1>0$, we put
\begin{equation}\label{e3.70}
W_1(T_1)=\left\{v\in L^{\infty}(0,T_1;H^1):v_{t}\in
L^{\infty}(0,T_1;L^2)\right\}.
\end{equation}
Then $W_1(T_1)$ is a Banach space with respect to the norm
(see[5]):
\begin{equation}\label{e3.71}
\|v\|_{W_1(T_1)}=\|v\|_{L^{\infty}(0,T_1;H^1)}+\|v_{t}\|_{L^{\infty}(0,T_1;L^2)},
v \in W_1(T_1).
\end{equation}
For $M>0$ and $T_1>0$, we put
\begin{equation}\label{e3.72}
B(M,T_1)=\left\{v\in W_1(T_1):\|v\|_{W_1(T_1)}\leq M\right\}.
\end{equation}
We also define the operator $\digamma$ from $B(M,T_1)$ into
$W(T_1)$ by $u=\digamma(v)$, where $u$ is the unique solution of
problem \eqref{e3.3}- \eqref{e3.5}. We would like to show that
$\digamma$ is a contraction operator from $B(M,T_1)$ into itself.
Applying the contraction mapping theorem, the operator $\digamma$
has a fixed point in $B(M,T_1)$ that is also a weak solution of
the problem \eqref{e1.1}- \eqref{e1.3}.

First, by Theorem 3.1, we note that the unique solution of problem
\eqref{e3.3}- \eqref{e3.5} satisfies \eqref{e3.7}, \eqref{e3.8},
\eqref{e3.9}. On the other hand, it follows from $(H3)$, that
\begin{equation}\label{e3.73}
\begin{aligned}
\int_0^{t}\|f(\cdot,s,v(s))\|^2ds &\leq 2\int_0^{t}\|f(\cdot,s,v(s))-f(\cdot,s,0)\|^2ds\\
&+2\int_0^{t}\|f(\cdot,s,0)\|^2ds\\
&\leq 2T_1K_1^2M^2+2\int_0^{T}\|f(\cdot,s,0)\|^2ds,
\end{aligned}
\end{equation}
where
\begin{equation}\label{e3.74}
\begin{aligned}
K_1&=K_1(M,T,f)\\
&=sup\left\{\left|D_3f(x,t,u)\right|:0\leq x\leq 1, 0\leq t\leq T,
|u|\leq \sqrt{2M}\right\}.
\end{aligned}
\end{equation}
It follows from \eqref{e3.7}-\eqref{e3.10} and \eqref{e3.73} that
\begin{equation}\label{e3.75}
\begin{aligned}
&\|u^{/}(t)\|^2+\|u(t)\|_{\eta}^2+2\int_0^{t}\|u^{/}(s)\|_{L^{q}}^{q}ds\\
&\leq \left(C_{1T}+2T_1K_1^2M^2\right)exp(T_1C_{2T}), \forall t\in
[0,T_1],
\end{aligned}
\end{equation}
where
\begin{equation}\label{e3.76}
\begin{gathered}
C_{1T}=C_{1T}\left(\widetilde{u}_0,\widetilde{u}_1,k,g\right)=2\left[
\|\widetilde{u}_1\|^2+\|\widetilde{u}_0\|_{\eta}^2+2\left|g_1(0)\widetilde{u}_0(1)\right|
\begin{array}{c}
   \\
   \\
\end{array}
\right.\\
\left.+6\|g_1\|_{L^{\infty}(0,T)}^2+2\|g_1^{/}\|_{L^2(0,T)}^2+2\int_0^{T}\|f(\cdot,s,0)\|^2ds\right],\\
C_{2T}=C_{2T}(k)=2\left[3+2|k(0)|+6\|k\|_{L^2(0,T)}^2+T\|k^{/}\|_{L^2(0,T)}^2\right].
\end{gathered}
\end{equation}
By choosing $M>0$ large enough so that $C_{1_T}=\frac{1}{4}M^2$,
then $T_1$ sufficiently small so that
\begin{equation}\label{e3.77}
\left(\frac{1}{4}M^2+2T_1K_1^2M^2\right)exp(T_1C_{2T})\leq
\frac{1}{2}M^2,
\end{equation}
and
\begin{equation}\label{e3.78}
2\sqrt{2T_1}K_1exp\left[T_1\left(2+2|k(0)|+2\|k\|_{L^2(0,T)}^2+\|k^{/}\|_{L^1(0,T)}^2\right)\right]<1.
\end{equation}
From \eqref{e3.75}, \eqref{e3.77} we have $\|u\|_{W_1(T_1)}\leq
M$, hence $u\in B(M,T_1)$. This shows that $\digamma$ maps
$B(M,T_1)$ into itself.

Next, we verify that $\digamma$ is a contraction. Let
$u_1=\digamma(v_1)$, $u_2=\digamma(v_2)$, where $v_1, v_2\in
B(M,T_1)$. Put $U=u_1- u_2$ and $V=v_1-v_2$. Then $U$ is the weak
solution of the following problem
\begin{equation}\label{e3.79}
\begin{aligned}
&U^{//}-U_{xx}+\int_0^{t}k(t-s)U_{xx}(s)ds+\psi_{q}(u_1^{/})-\psi_{q}(u_2^{/})\\
&=f(x,t,v_1(t))-f(x,t,v_2(t)), 0<x<1, 0<t<T_1,\\
&U_{x}(0,t)-U(0,t)=U_{x}(1,t)+\eta U(1,t)=0,\\
&U(0)=U^{/}(0)=0, \\
&U\in W(T_1);  U^{//},  U_{xx}\in L^{q^{/}}(0,T_1;(H^1)^{/}).
\end{aligned}
\end{equation}
By using Lemma 3.2 with $\widetilde{u}_0=\widetilde{u}_1=0$,
$\Phi=-\psi_{q}(u_1^{/})+\psi_{q}(u_2^{/})+f(x,t,v_1(t))-f(x,t,v_2(t))$,
we have
\begin{equation}\label{e3.80}
\begin{aligned}
\delta(t)=&-2k(0)\int_0^{t}\|U(r)\|_{\eta}^2dr+2\int_0^{t}k(t-s)a(U(s),U(t))ds\\
&-2 \int_0^{t}dr\int_0^{r}k^{/}(r-s)a(U(s),U(r))ds\\
&+2 \int_0^{t}<f(\cdot,s,v_1(s))-f(\cdot,s,v_2(s)),U^{/}(s)>ds,
 \text{a.e.} \quad t\in [0,T1],
\end{aligned}
\end{equation}
where
\begin{equation}\label{e3.81}
\begin{aligned}
\delta(t)&=\|U^{/}(t)\|^2+\|U(t)\|_{\eta}^2+2\int_0^{t}\left<\psi_{q}(u_1^{/})-\psi_{q}(u_2^{/}),U^{/}(s)
\right>ds\\
&\geq \|U^{/}(t)\|^2+\|U(t)\|_{\eta}^2+2C_{q}
\int_0^{t}\|U^{/}(s)\|_{L^{q}}^{q}ds.
\end{aligned}
\end{equation}
By the assumption $(H4)$, we have
\begin{equation}\label{e3.82}
\begin{aligned}
&2\int_0^{t}\left<f(\cdot,s,v_1(s))-f(\cdot,s,v_2(s)),U^{/}(s)\right>ds\\
&\leq
\int_0^{t}\|U^{/}(s)\|^2ds+\int_0^{t}\|f(\cdot,s,v_1(s))-f(\cdot,s,v_2(s))\|^2ds\\
&\leq\int_0^{t}\|U^{/}(s)\|^2ds+2T_1K_1^2\|V\|_{W_1(T_1)}^2,
\end{aligned}
\end{equation}
Therefore, we can prove in a similar manner as above that
\begin{equation}\label{e3.83}
\begin{aligned}
\delta(t)\leq & 2T_1K_1^2\|V\|_{W_1(T_1)}^2\\
&+2\left(2+2|k(0)|+2\|k\|_{L^2(0,T)}^2+\|k^{/}\|_{L^1(0,T)}^2\right)\int_0^{t}\delta(s)ds.
\end{aligned}
\end{equation}
By Gronwall's lemma, we obtain from \eqref{e3.83} that
\begin{equation}\label{e3.84}
\delta(t)=2\left(\rho_1(k,K_1,T,T_1)\|V\|_{W_1(T_1)}\right)^2,
\end{equation}
where
\begin{equation}\label{e3.85}
\begin{aligned}
\rho_1&(k,K_1,T,T_1)=\\
&\sqrt{2T_1}K_1exp\left[T_1\left(2+2|k(0)|+2\|k\|_{L^2(0,T)}^2+\|k^{/}\|_{L^1(0,T)}^2\right)\right].
\end{aligned}
\end{equation}
It follows from \eqref{e3.81}, \eqref{e3.84} and \eqref{e3.85}
that
\begin{equation}\label{e3.86}
\|U\|_{W_1(T_1)}\leq 2\rho_1(k,K_1,T,T_1)\|V\|_{W_1(T_1)},
\end{equation}
where
\begin{equation}\label{e3.87}
2\rho_1(k,K_1,T,T_1)<1,
\end{equation}
since \eqref{e3.78} and \eqref{e3.85}.

Hence, \eqref{e3.86} shows that $\digamma:B(M,T_1)\rightarrow
B(M,T_1)$ is a contraction. Applying the contraction mapping
theorem, the operator $\digamma$ has a fixed point in $B(M,T_1)$
that is also a weak solution of the problem \eqref{e1.1}-
\eqref{e1.3}.

The solution of the problem \eqref{e1.1}- \eqref{e1.3} is unique,
that can be showed using the same arguments as in the proof of
Theorem 3.1. The proof of Theorem 3.3 is completed.
\end{proof}
\begin{remark}\label{rmk3.4}
 In the case of $\lambda =0$, $f(x,t,u)=|u|^{p-2}u$, $p>2$, $k\in W^{2,1}\left(\mathbb{R}_+\right)$,
 $k\geq 0$, $k(0)>0$, $0<\int_0^{+\infty}k(t)dt<1$, $k^{/}(t)+\zeta k(t)\leq 0$ for all $t\geq 0$,
 with $\zeta>0$, and the boundary condition $u(0,t)=u(1,t)=0$ standing for \eqref{e1.2}, S. Berrimia,
 S. A. Messaoudi [1] has obtained a global existence and uniqueness theorem.
\end{remark}

\section{Decay of solution}
In this part, we will consider the problem of global existence and
asymptotic behavior for $t\rightarrow +\infty$. We assume that
$g(t)=0$, $f(x,t,u)=F(x,t) -|u|^{p-2}u$, $p\geq 2$ and consider
the following problem
\begin{equation}\label{e4.1}
\begin{aligned}
&u_{tt}-u_{xx}+\int_0^{t}k(t-s)u_{xx}(s)ds+|u|^{p-2}u+|u_{t}|^{q-2}u_{t}\\
&=F(x,t), 0<x<1, t>0,\\
&u_{x}(0,t)=u(0,t),  u_{x}(1,t)+\eta u(1,t)=0,\\
&u(x,0)=\widetilde{u}_0(x),  u_{t}(x,0)=\widetilde{u}_1(x),
\end{aligned}
\end{equation}
We make the following assumptions:
\begin{itemize}
\item[($\widetilde{H}1$)] $\eta \geq 0, p, q\geq 2$,
\item[($\widetilde{H}2)$] $k \in W^{2,1}(\mathbb{R}_+)$, $k\geq 0$, satisfying
  \begin{itemize}
  \item[(i)] $k(0)>0$, $0<1-\int_0^{+\infty}k(t)dt=k_{\infty}<1$,
  \item[(ii)] there exists a positive constant $\zeta$ such that $k^{/}(t)+\zeta k(t)\leq 0$ for all $t\geq 0$,
  \end{itemize}
\item[($\widetilde{H}3)$] $\widetilde{u}_0\in H^2$ and $\widetilde{u}_1\in H^1$,
\item[($\widetilde{H}4)$] $F\in L^1(0,\infty;L^2)\cap L^2(0,\infty;L^2)$, $F_{t}\in L^1(0,\infty;L^2)$,
\item[($\widetilde{H}5)$] There exists a constant $\sigma>0$ such that $\int_0^{\infty}e^{st}\|F(t)\|^2dt<+\infty$.
\end{itemize}
Under assumptions $(\widetilde{H}1)$-$(\widetilde{H}4)$ and let
$T>0$, by theorem 2.3, the problem \eqref{e4.1} has a unique weak
solution $u(t)$ such that
\begin{equation}\label{e4.2}
u\in L^{\infty}(0,T;H^2), u_{t}\in L^{\infty}(0,T;H^1), u_{tt}\in
L^{\infty}(0,T;L^2).
\end{equation}
Then, we have the following
\begin{lemma}\label{lem4.1}
Suppose that $(\widetilde{H}1)-(\widetilde{H}4)$ hold. Then there
is a unique solution $u(t)$ of problem \eqref{e4.1} defined on
$\mathbb{R}_+$. Moreover
\begin{equation}\label{e4.3}
\|u^{/}(t)\| +\|u(t)\|_{\eta}\leq C \quad \text{for all} \quad
t\geq 0,
\end{equation}
where $C$ is a positive constant depending only on
$\widetilde{u}_0$, $\widetilde{u}_1$, $F$, $k_{\infty}$ and $p$.
\end{lemma}
\begin{proof}
By multiplying the equation \eqref{e4.1}$_1$ by $u_{t}$ and
integrate over $(0,1)\times (0,t)$ we obtain
\begin{equation}\label{e4.4}
\begin{aligned}
E(t)&+2\int_0^{t}\|u^{/}(s)\|_{L^{q}}^{q}ds+\int_0^{t}k(s)\|u(s)\|_{\eta}^2ds\\
&-\int_0^{t}dr\int_0^{r}k^{/}(r-s)\|u(s)-u(r)\|_{\eta}^2ds\\
&=E(0)+2\int_0^{t}<F(s),u^{/}(s)>ds,
\end{aligned}
\end{equation}
where
\begin{equation}\label{e4.5}
\begin{aligned}
E(t)=&\|u^{/}(t)\|^2+\left(1-\int_0^{t}k(s)ds\right)\|u(t)\|_{\eta}^2+\frac{2}{p}\|u(t)\|_{L^{p}}^{p}\\
&+ \int_0^{t}k(t-s)\|u(s)-u(t)\|_{\eta}^2ds.
\end{aligned}
\end{equation}
On the other hand, by $(\widetilde{H}4)$ and the Cauchy's
inequality, we obtain
\begin{equation}\label{e4.6}
\begin{aligned}
2\int_0^{t}\left<F(s),u^{/}(s)\right>ds&\leq
\int_0^{t}\|F(s)\|ds+\int_0^{t}\|F(s)\|\|u^{/}(s)\|^2ds\\
& \leq \int_0^{+\infty}\|F(s)\|ds+\int_0^{t}\|F(s)\|E(s)ds.
\end{aligned}
\end{equation}
By Gronwall's lemma, we obtain from \eqref{e4.4} and \eqref{e4.6}
that
\begin{equation}\label{e4.7}
\begin{aligned}
E(t)& \leq
\left(E(0)+\int_0^{+\infty}\|F(s)\|ds\right)exp\left(\int_0^{t}\|F(s)\|ds\right)\\
&\leq
\left(E(0)+\int_0^{+\infty}\|F(s)\|ds\right)exp\left(\int_0^{+\infty}\|F(s)\|ds\right)=C,
\forall t\geq 0.
\end{aligned}
\end{equation}
By $(\widetilde{H}3,i)$, we have
\begin{equation}\label{e4.8}
E(t)\geq
\|u^{/}(t)\|^2+\left(1-\int_0^{t}k(s)ds\right)\|u(t)\|_{\eta}^2\geq
\|u^{/}(t)\|^2+k_{\infty}\|u(t)\|_{\eta}^2.
\end{equation}
Then we obtain \eqref{e4.3} from \eqref{e4.7} and \eqref{e4.8}.
This completes the proof of Lemma 4.1.
\end{proof}
In this section we state and prove decay result.
\begin{theorem}\label{thm4.2}
Suppose that $(\widetilde{H}1)-(\widetilde{H}5)$ hold. Then the
solution $u(t)$ of problem \eqref{e4.1} decays exponentially to
zero as $t\rightarrow +\infty$ in the following sense: there exist
the positive constants $N$ and $\gamma$ such that
\begin{equation}\label{e4.9}
\|u^{/}(t)\| +\|u(t)\|_{\eta}\leq Ne^{-\gamma t} \quad \text{for
all} \quad  t\geq 0.
\end{equation}
\end{theorem}
\begin{proof}
We use the following functional
\begin{equation}\label{e4.10}
\Gamma(t)=\Gamma(\varepsilon_1,\varepsilon_2,t)=E(t)+\varepsilon_1E_1(t)+\varepsilon_2E_2(t),
\end{equation}
where
\begin{gather}\label{e4.11}
E_1(t)=<u(t),u^{/}(t)>,\\
\label{e4.12}
E_2(t)=-\int_0^{t}k(t-s)\left<u^{/}(t),u(t)-u(s)\right>ds.
\end{gather}
{\emph{Estimating}} $\Gamma(t)$. \\
By \eqref{e2.3}, \eqref{e2.4}, we obtain from $(\widetilde{H}2,
i)$ that
\begin{equation}\label{e4.13}
|E_1(t)|=\left|<u(t),u^{/}(t)>\right|\leq\frac{1}{2}\|u^{/}(t)\|^2+\|u(t)\|_{\eta}^2,
\end{equation}
\begin{equation}\label{e4.14}
\begin{aligned}
|E_2(t)|&=\left|\int_0^{t}k(t-s)\left<u^{/}(t),u(t)-u(s)\right>ds\right|\\
&\leq
\frac{1}{2}\|u^{/}(t)\|^2+\frac{1}{2}\left(\int_0^{t}k(t-s)\|u(t)-u(s)\|ds\right)^2\\
&\leq
\frac{1}{2}\|u^{/}(t)\|^2+(1-k_{\infty})\int_0^{t}k(t-s)\|u(t)-u(s)\|_{\eta}^2ds.
\end{aligned}
\end{equation}
Hence, it follows from \eqref{e4.10}-\eqref{e4.14} that for
$\varepsilon_1$, $\varepsilon_2$ small enough, there exist two
positive constants $\alpha_1, \alpha_2$, such that
\begin{equation}\label{e4.15}
\alpha_1E(t)\leq \Gamma(t)\leq \alpha_2E(t).
\end{equation}
{\emph{Estimating} } $\Gamma^{/}(t)$. \\
Now differentiating \eqref{e4.4} with respect to $t$, we have
\begin{equation}\label{e4.16}
\begin{aligned}
E^{/}(t)&=-2\|u^{/}(t)\|_{L^{q}}^{q}+\int_0^{t}k^{/}(t-s)\|u(s)-u(t)\|_{\eta}^2ds\\
&-k(t)\|u(t)\|_{\eta}^2+2\left<F(t),u^{/}(t)\right>\\
&\leq
-2\|u^{/}(t)\|_{L^{q}}^{q}+\int_0^{t}k^{/}(t-s)\|u(s)-u(t)\|_{\eta}^2ds+2\left<F(t),u^{/}(t)\right>,
\end{aligned}
\end{equation}
since  $k(t)\geq 0$.

By multiplying the equation \eqref{e4.1}$_1$ by $u$ and integrate
over $(0,1)$ we obtain
\begin{equation}\label{e4.17}
\begin{aligned}
E_1^{/}(t)&=\|u^{/}(t)\|^2-\|u(t)\|_{\eta}^2-\|u(t)\|_{L^{p}}^{p}+\left<F(t),u(t)\right>\\
&+\int_0^{t}k(t-s)a(u(s),u(t))ds-\left<|u^{/}(t)|^{q-2}u^{/}(t),u(t)\right>\\
&=\|u^{/}(t)\|^2-\|u(t)\|_{\eta}^2-\|u(t)\|_{L^{p}}^{p}+\left<F(t),u(t)\right>+I_1(t)+I_2(t).
\end{aligned}
\end{equation}
We now estimate the last two terms in the right side of
\eqref{e4.17} as follows
\vspace{0.2cm}\\
{\emph{Estimating}} $I_1(t)$.\\
Using the inequality
\begin{equation}\label{e4.18}
ab\leq
\frac{\delta}{r}a^{r}+\frac{r-1}{r}\delta^{\frac{-r}{r-1}}b^{\frac{r}{r-1}},
\forall a,b\geq 0, \forall r>1, \forall \delta>0,
\end{equation}
we have
\begin{equation}\label{e4.19}
\begin{aligned}
I_1(t)&=\int_0^{t}k(t-s)a(u(s),u(t))ds\\
&=\int_0^{t}k(t-s)a\left(u(s)-u(t),u(t)\right)ds+\int_0^{t}k(t-s)\|u(t)\|_{\eta}^2ds\\
&\leq
\delta_1\|u(t)\|_{\eta}^2+\frac{1}{4\delta_1}\left(\int_0^{t}k(s)ds\right)\left(\int_0^{t}k(t-s)\|u(s)-u(t)\|_{\eta}^2ds\right)\\
& +\left(\int_0^{t}k(s)ds\right)\|u(t)\|_{\eta}^2\\
&\leq
\delta_1\|u(t)\|_{\eta}^2+\frac{1-k_\infty}{4\delta_1}\int_0^{t}k(t-s)\|u(s)-u(t)\|_{\eta}^2ds\\
&+(1-k_{\infty})\|u(t)\|_{\eta}^2\\
&\leq
(\delta_1+1-k_{\infty})\|u(t)\|_{\eta}^2+\frac{1-k_\infty}{4\delta_1}\int_0^{t}k(t-s)\|u(s)-u(t)\|_{\eta}^2ds,
\end{aligned}
\end{equation}
for all $\delta_1>0$.\\
{\emph{Estimating}} $I_2(t)$.\\
We again use inequality \eqref{e4.18} we obtain from \eqref{e4.3}
that
\begin{equation}\label{e4.20}
\begin{aligned}
I_2(t)&=-\left<|u^{/}(t)|^{q-2}u^{/}(t),u(t)\right>\leq\|u^{/}(t)\|_{L^{q}}^{q-1}\|u(t)\|_{L^{q}}\\
&\leq\frac{\delta_1^q}{q}\|u(t)\|_{L^{q}}^{q}+\frac{q-1}{q}\delta_1^{\frac{-q}{q-1}}\|u^{/}(t)\|_{L^{q}}^{q}\\
&\leq
2\frac{\delta_1^q}{q}\left(\sqrt{2}C\right)^{q-2}\|u(t)\|_{\eta}^2+\frac{q-1}{q}\delta_1^{\frac{-q}{q-1}}\|u^{/}(t)\|_{L^{q}}^{q},
\end{aligned}
\end{equation}
for all $\delta_1>0$.\\
By combining \eqref{e4.17}, \eqref{e4.19} and \eqref{e4.20}, we
obtain
\begin{equation}\label{e4.21}
\begin{aligned}
E_1^{/}(t)&\leq -\|u(t)\|_{L^{p}}^{p}+\|u^{/}(t)\|^2
-\left(k_{\infty}-\delta_1-2\frac{\delta_1^q}{q}\left(\sqrt{2}C\right)^{q-2}\right)\|u(t)\|_{\eta}^2\\
&+\frac{q-1}{q}\delta_1^{\frac{-q}{q-1}}\|u^{/}(t)\|_{L^{q}}^{q}
+\frac{1-k_\infty}{4\delta_1}\int_0^{t}k(t-s)\|u(s)-u(t)\|_{\eta}^2ds\\
&+\left<F(t),u(t)\right>.
\end{aligned}
\end{equation}
Then, we can always choose the constant $\delta_1>0$ such that
\begin{equation}\label{e4.22}
\gamma_1=k_{\infty}-\delta_1-2\frac{\delta_1^q}{q}\left(\sqrt{2}C\right)^{q-2}>0.
\end{equation}
This implies that
\begin{equation}\label{e4.23}
\begin{aligned}
E_1^{/}(t)&\leq
-\|u(t)\|_{L^{p}}^{p}+\|u^{/}(t)\|^2-\gamma_1\|u(t)\|_{\eta}^2+\gamma_2\|u^{/}(t)\|_{L^{q}}^{q}\\
&+\gamma_3\int_0^{t}k(t-s)\|u(s)-u(t)\|_{\eta}^2ds+\left<F(t),u(t)\right>,
\end{aligned}
\end{equation}
where
\begin{equation}\label{e4.24}
\gamma_2=\frac{q-1}{q}\delta_1^{\frac{-q}{q-1}},
\gamma_3=\frac{1-k_\infty}{4\delta_1}.
\end{equation}
Direct calculations give
\begin{equation}\label{e4.25}
\begin{aligned}
E_2^{/}(t)&=-\left(\int_0^{t}k(s)ds\right)\|u^{/}(t)\|^2-\int_0^{t}k^{/}(t-s)\left<u^{/}(t),u(t)-u(s)\right>ds\\
&+\int_0^{t}k(t-s)a(u(t),u(t)-u(s))ds\\
&-\int_0^{t}k(t-s)a\left(\int_0^{t}k(t-\tau)u(\tau)d\tau,u(t)-u(s)\right)ds\\
&+\int_0^{t}k(t-s)\left<|u(t)|^{p-2}u(t),u(t)-u(s)\right>ds\\
&+\int_0^{t}k(t-s)\left<|u^{/}(t)|^{q-2}u^{/}(t),u(t)-u(s)\right>ds\\
&-\int_0^{t}k(t-s)\left<F(t),u(t)-u(s)\right>ds=\sum_{i=1}^{7}J_{i}(t).
\end{aligned}
\end{equation}
Similarly to \eqref{e4.17}, we estimate respectively the following
terms on the right-hand side of \eqref{e4.25} as follows.\\
{\emph{Estimating}} $J_1(t)$.\\
Since $k$ is continuous and $k(0)>0$ then there exists $t_0>0$,
such that
\begin{equation}\label{e4.26}
\int_0^{t}k(s)ds\geq \int_0^{t_0}k(s)ds=k_0>0 \quad \text{for all}
\quad t\geq t_0.
\end{equation}
Hence,
\begin{equation}\label{e4.27}
J_1(t)=-\left(\int_0^{t}k(s)ds\right)\|u^{/}(t)\|^2\leq
-k_0\|u^{/}(t)\|^2 \quad \text{for all} \quad t\geq t_0.
\end{equation}
{\emph{Estimating}} $J_2(t)$.\\
\begin{equation}\label{e4.28}
\begin{aligned}
&J_2(t)=-\int_0^{t}k^{/}(t-s) \left<u^{/}(t),u(t)-u(s)\right>ds\\
&\leq
\delta_2\|u^{/}(t)\|^2+\frac{1}{4\delta_2}\left(\int_0^{t}|k^{/}(t-s)|ds\right)
\left(\int_0^{t}|k^{/}(t-s)|\|u(s)-u(t)\|^2ds\right)\\
&\leq
\delta_2\|u^{/}(t)\|^2+\frac{1}{2\delta_2}\left(\int_0^{t}|k^{/}(t-s)|ds\right)
\left(\int_0^{t}|k^{/}(t-s)|\|u(s)-u(t)\|_{\eta}^2ds\right)\\
&\leq
\delta_2\|u^{/}(t)\|^2-\frac{k(0)}{2\delta_2}\int_0^{t}k^{/}(t-s)\|u(s)-u(t)\|_{\eta}^2ds.\\
\end{aligned}
\end{equation}
{\emph{Estimating}} $J_3(t)$.\\
\begin{equation}\label{e4.29}
\begin{aligned}
J_3(t)&=\int_0^{t}k(t-s)a\left(u(t),u(t)-u(s)\right)ds\\
&\leq
\delta_2\|u(t)\|_{\eta}^2+\frac{1}{4\delta_2}\left(\int_0^{t}k(s)ds\right)
\left(\int_0^{t}k(t-s)\|u(s)-u(t)\|_{\eta}^2ds\right)\\
&\leq
\delta_2\|u(t)\|_{\eta}^2+\frac{1-k_\infty}{4\delta_2}\int_0^{t}k(t-s)\|u(s)-u(t)\|_{\eta}^2ds
\end{aligned}
\end{equation}
{\emph{Estimating}} $J_4(t)$.\\
\begin{equation}\label{e4.30}
\begin{aligned}
J_4(t)&=-\int_0^{t}k(t-s)a\left(\int_0^{t}k(t-\tau)u(\tau)d\tau,u(t)-u(s)\right)ds\\
&\leq \int_0^{t}k(t-\tau)\|u(\tau)\|_{\eta}d\tau
\int_0^{t}k(t-s)\|u(s)-u(t)\|_{\eta}ds\\
&\leq
\delta_2\left(\int_0^{t}k(t-\tau)\|u(\tau)\|_{\eta}d\tau\right)^2\\
&+\frac{1}{4\delta_2}\left(\int_0^{t}k(t-s)\|u(s)-u(t)\|_{\eta}ds\right)^2\\
&\leq
2\delta_2\left(\int_0^{t}k(t-\tau)\|u(\tau)\|_{\eta}d\tau\right)^2\\
&+\left(2\delta_2+\frac{1}{4\delta_2}\right)\left(\int_0^{t}k(t-s)\|u(s)-u(t)\|_{\eta}ds\right)^2\\
&\leq 2\delta_2(1-k_{\infty})^2\|u(t)\|_{\eta}^2\\
&+\left(2\delta_2+\frac{1}{4\delta_2}\right)(1-k_{\infty})
\int_0^{t}k(t-s)\|u(s)-u(t)\|_{\eta}^2ds.
\end{aligned}
\end{equation}
{\emph{Estimating}} $J_5(t)$.
\begin{equation}\label{e4.31}
\begin{aligned}
&J_5(t)=\int_0^{t}k(t-s)\left<|u(t)|^{p-2}u(t),u(t)-u(s)\right>ds\\
&\leq
2\left(\sqrt{2}C\right)^{p-2}\int_0^{t}k(t-s)\|u(t)\|_{\eta}\|u(t)-u(s)\|_{\eta}ds\\
&\leq 2\left(\sqrt{2}C\right)^{p-2}\left[
\delta_2\|u(t)\|_{\eta}^2+\frac{1}{4\delta_2}\left(\int_0^{t}k(t-s)
\|u(t)-u(s)\|_{\eta}ds\right)^2\right]\\
&\leq 2\left(\sqrt{2}C\right)^{p-2}\left[
\delta_2\|u(t)\|_{\eta}^2+\frac{1}{4\delta_2}(1-k_{\infty})\int_0^{t}k(t-s)
\|u(t)-u(s)\|_{\eta}^2ds\right].
\end{aligned}
\end{equation}
{\emph{Estimating}} $J_6(t)$. We again use inequality
\eqref{e4.18} with $r=q$, $\delta= \delta_2$, we obtain from
\eqref{e4.3} that
\begin{equation}\label{e4.32}
\begin{aligned}
&\left<|u^{/}(t)|^{q-2}u^{/}(t),u(t)-u(s)\right>\leq\|u^{/}(t)\|_{L^{q}}^{q-1}\|u(t)-u(s)\|_{L^{q}}\\
&\leq\frac{\delta_2^{q}}{q}\|u(t)-u(s)\|_{L^{q}}^{q}+
\frac{q-1}{q}\delta_2^{\frac{-q}{q-1}}\|u^{/}(t)\|_{L^{q}}^{q}\\
&\leq
2\frac{\delta_2^{q}}{q}\left(2\sqrt{2}C\right)^{q-2}\|u(t)-u(s)\|_{\eta}^2
+\frac{q-1}{q}\delta_2^{\frac{-q}{q-1}}\|u^{/}(t)\|_{L^{q}}^{q}.
\end{aligned}
\end{equation}
It follows from \eqref{e4.32} that
\begin{equation}\label{e4.33}
\begin{aligned}
J_6(t)&=\int_0^{t}k(t-s)\left<|u^{/}(t)|^{q-2}u^{/}(t),u(t)-u(s)\right>ds\\
&\leq
2\frac{\delta_2^{q}}{q}\left(2\sqrt{2}C\right)^{q-2}\int_0^{t}k(t-s)
\|u(t)-u(s)\|_{\eta}^2ds\\
&+\frac{q-1}{q}\delta_2^{\frac{-q}{q-1}}\|u^{/}(t)\|_{L^{q}}^{q}\int_0^{t}k(t-s)
ds\\
&\leq
2\frac{\delta_2^{q}}{q}\left(2\sqrt{2}C\right)^{q-2}\int_0^{t}k(t-s)
\|u(t)-u(s)\|_{\eta}^2ds\\
&+\frac{q-1}{q}\delta_2^{\frac{-q}{q-1}}(1-k_\infty)\|u^{/}(t)\|_{L^{q}}^{q}
\end{aligned}
\end{equation}
{\emph{Estimating}} $J_7(t)$
\begin{equation}\label{e4.34}
\begin{aligned}
J_7(t)&=-\int_0^{t}k(t-s)\left<F(t),u(t)-u(s)\right>ds\\
&\leq \int_0^{t}k(t-s)\|F(t)\|\|u(t)-u(s)\|ds\\
&\leq
\frac{1}{4\delta_2}\|F(t)\|^2+\delta_2\left(\int_0^{t}k(t-s)ds\right)
\left(\int_0^{t}k(t-s)\|u(t)-u(s)\|^2ds\right)\\
&\leq \frac{1}{4\delta_2}\|F(t)\|^2
+2\delta_2(1-k_{\infty})\int_0^{t}k(t-s)\|u(t)-u(s)\|_{\eta}^2ds.
\end{aligned}
\end{equation}
By combining \eqref{e4.25}, \eqref{e4.27}-\eqref{e4.31},
\eqref{e4.33} and \eqref{e4.34}, we obtain
\begin{equation}\label{e4.35}
\begin{aligned}
E_2^{/}(t)&=-\left(k_0-\delta_2\right)\|u^{/}(t)\|^2+\delta_2\widehat{\gamma}_1\|u(t)\|_{\eta}^2
+\widehat{\gamma}_2\|u^{/}(t)\|_{L^{q}}^{q}\\
&+\widehat{\gamma}_3\int_0^{t}k(t-s)\|u(t)-u(s)\|_{\eta}^2ds\\
&-\widehat{\gamma}_4\int_0^{t}k^{/}(t-s)\|u(s)-u(t)\|_{\eta}^2ds
+\frac{1}{4\delta_2}\|F(t)\|^2,
\end{aligned}
\end{equation}
where
\begin{equation}\label{e4.36}
\begin{gathered}
\widehat{\gamma}_1=1+2(1-k_{\infty})^2+2\left(\sqrt{2}C\right)^{p-2},\\
\widehat{\gamma}_2=\frac{q}{q-1}\delta_2^{\frac{-q}{q-1}}(1-k_{\infty}),\\
\widehat{\gamma}_3=2\frac{\delta_2^q}{q}\left(2\sqrt{2}C\right)^{q-2}
+(1-k_{\infty})\left[\frac{1}{2\delta_2}\left(\sqrt{2}C\right)^{p-2}+\left(4\delta_2
+\frac{1}{2\delta_2}\right)\right],\\
\widehat{\gamma}_4=\frac{k(0)}{2\delta_2}.
\end{gathered}
\end{equation}
Combining of \eqref{e4.10}, \eqref{e4.16}, \eqref{e4.23} and
\eqref{e4.35}, we obtain
\begin{equation}\label{e4.37}
\begin{aligned}
\Gamma^{/}(t)&+\varepsilon_1\|u(t)\|_{L^{p}}^{p}+\left((k_0-\delta_2)\varepsilon_2-\varepsilon_1\right)
\|u^{/}(t)\|^2\\
&+\left(\varepsilon_1\gamma_1-\varepsilon_2\delta_2\widehat{\gamma}_1\right)
\|u(t)\|_{\eta}^2+\left(2-\varepsilon_1\gamma_2-\varepsilon_2\widehat{\gamma}_2\right)
\|u^{/}(t)\|_{L^{q}}^{q}\\
&-\left(\varepsilon_1\gamma_3+\varepsilon_2\widehat{\gamma}_3\right)
\int_0^{t}k(t-s)\|u(t)-u(s)\|_{\eta}^2ds\\
&-\left(1-\varepsilon_2\widehat{\gamma}_4\right)\int_0^{t}k^{/}(t-s)\|u(s)-u(t)\|_{\eta}^2ds\\
&\leq
\left<F(t),2u^{/}(t)+\varepsilon_1u(t)\right>+\frac{\varepsilon_2}{4\delta_2}\|F(t)\|^2.
\end{aligned}
\end{equation}
Whence $\delta_1$ is fixed, choosing
\begin{equation}\label{e4.38}
\delta_2=\frac{1}{2}\frac{k(0)\gamma_1}{\gamma_1+\widehat{\gamma}_1},
\varepsilon_2=\frac{2}{k_0}\varepsilon_1, \quad \text{where} \quad
\varepsilon_1>0 \quad \text{is arbitrary},
\end{equation}
we deduce from \eqref{e4.37} and \eqref{e4.38} that
\begin{equation}\label{e4.39}
\begin{aligned}
\Gamma^{/}(t)&+\varepsilon_1\|u(t)\|_{L^{p}}^{p}+\frac{\varepsilon_1\widehat{\gamma}_1}
{\gamma_1+\widehat{\gamma}_1}\|u^{/}(t)\|^2+\frac{\varepsilon_1\gamma_1^2}
{\gamma_1+\widehat{\gamma}_1}\|u(t)\|_{\eta}^2\\
&+\left(2-\varepsilon_1\left(1+\frac{2}{k(0)}\widehat{\gamma}_2\right)
\right)\|u^{/}(t)\|_{L^{q}}^{q}\\
&-\varepsilon_1\left(\gamma_3+\frac{2}{k(0)}\widehat{\gamma}_3\right)
\int_0^{t}k(t-s)\|u(t)-u(s)\|_{\eta}^2ds\\
&-\left(1-\frac{2}{k(0)}\varepsilon_1\widehat{\gamma}_4\right)\int_0^{t}k^{/}(t-s)\|u(s)-u(t)\|_{\eta}^2ds\\
&\leq
\left<F(t),2u^{/}(t)+\varepsilon_1u(t)\right>+\frac{\varepsilon_1}{k_0^2}
\left(1+\frac{\widehat{\gamma}_1}{\gamma_1}\right)\|F(t)\|^2.
\end{aligned}
\end{equation}
Next, we choose $\varepsilon_1>0$, with
$$
\varepsilon_1<min\left\{\frac{\zeta}{\gamma_3+\frac{2}{k_0}\widehat{\gamma}_3
+\frac{2}{k_0}\widehat{\gamma}_4\zeta},
\frac{2}{1+\frac{2}{k_0}\widehat{\gamma}_2}\right\}
$$
and \eqref{e4.15} is satisfied, then by the assumption
$(\widetilde{H}2, ii)$, we deduce that
\begin{equation}\label{e4.40}
\begin{aligned}
\Gamma^{/}(t)+\varepsilon_1\|u(t)\|_{L^{p}}^{p}&+\frac{\varepsilon_1\widehat{\gamma}_1}
{\gamma_1+\widehat{\gamma}_1}\|u^{/}(t)\|^2+\frac{\varepsilon_1\gamma_1^2}
{\gamma_1+\widehat{\gamma}_1}\|u(t)\|_{\eta}^2+k_1\|u^{/}(t)\|_{L^{q}}^{q}\\
&+k_2\int_0^{t}k(t-s)\|u(s)-u(t)\|_{\eta}^2ds\\
&\leq \left<F(t),2u^{/}(t)+\varepsilon_1u(t)\right>+k_3\|F(t)\|^2,
\end{aligned}
\end{equation}
where
\begin{equation}\label{e4.41}
\begin{gathered}
k_1=2-\varepsilon_1\left(1+\frac{2}{k_0}\widehat{\gamma}_2\right) >0,\\
k_2=\zeta\left(1-\frac{2}{k_0}\varepsilon_1\widehat{\gamma}_4\right)
-\varepsilon_1\left(\gamma_3+\frac{2}{k_0}\widehat{\gamma}_3\right)>0,\\
k_3=\frac{\varepsilon_1}{k_0^2}\left(1+\frac{\widehat{\gamma}_1}{\gamma_1}\right).
\end{gathered}
\end{equation}
By combining \eqref{e4.5}, \eqref{e4.15} and \eqref{e4.40}, we can
always choose the constant $\widetilde{\gamma}>0$ is independent
of $t$ such that
\begin{equation}\label{e4.42}
\Gamma^{/}(t)+2\widetilde{\gamma}\Gamma(t)\leq
\left<F(t),2u^{/}(t)+\varepsilon_1u(t)\right>+k_3\|F(t)\|^2,
\end{equation}
for all $t\geq t_0$.\\
On ther other hand,
\begin{equation}\label{e4.43}
\left<F(t),2u^{/}(t)+\varepsilon_1u(t)\right>+k_3\|F(t)\|^2\leq
\widetilde{N}\|F(t)\|^2+\widetilde{\gamma}\Gamma(t),
\end{equation}
for some constant $\widetilde{N}>0$. Therefore
\begin{equation}\label{e4.44}
\Gamma^{/}(t)+\widetilde{\gamma}\Gamma(t)\leq
\widetilde{N}\|F(t)\|^2 \quad \text{for all} \quad t\geq t_0.
\end{equation}
Putting $\gamma=\frac{1}{2}min\{\sigma,\widetilde{\gamma}\}$. A
simple integration of \eqref{e4.44} over $(t_0, t)$ gives
\begin{equation}\label{e4.45}
\Gamma(t)\leq \left[\Gamma(t_0)e^{\sigma
t_0}+\widetilde{N}\int_{t_0}^{+\infty}e^{\sigma
s}\|F(s)\|^2ds\right]e^{-2\gamma t}=N_1e^{-2\gamma t},
\end{equation}
for all $t\geq t_0$.\\
By the boundedness of $\Gamma(t)$ on $[0, t_0]$, we deduce from
\eqref{e4.45} that
\begin{equation}\label{e4.46}
\Gamma(t)=\|\Gamma\|_{L^{\infty}(0,t_0)}e^{-2\gamma(t-t_0)}+N_1e^{-2\gamma
t}=N_2e^{-2\gamma t},
\end{equation}
for all $t\geq 0$.\\
By \eqref{e4.15}, it follows from \eqref{e4.46} that
\begin{equation}\label{e4.47}
E(t)\leq \frac{1}{\alpha_1}\Gamma(t)\leq
\frac{1}{\alpha_1}N_2e^{-2\gamma t}, \quad \text{for all} \quad
t\geq 0.
\end{equation}
This completes the proof of Theorem 4.2.
\end{proof}
\begin{remark}\label{rmk4.3}
The estimate \eqref{e4.9} holds for any regular solution
corresponding to $(\widetilde{u}_0,\widetilde{u}_1)\in H^2 \times
H^1$. This remains holds for solutions corresponding to
$(\widetilde{u}_0,\widetilde{u}_1)\in H^1 \times L^2$ by simple
density argument.
\end{remark}

\section{Numerical results}
Consider the following problem:
\begin{equation}\label{e5.1}
u_{tt}-u_{xx}+\int_0^{t}k(t-s)u_{xx}(s)ds+u_{t}^3=u^2+F(x,t),
0<x<1, 0<t<T,
\end{equation}
with boundary conditions
\begin{equation}\label{e5.2}
u_{x}(0,t)=u(0,t),  u_{x}(1,t)+u(1,t)=0,
\end{equation}
and initial conditions
\begin{equation}\label{e5.3}
u(x,0)=\widetilde{u}_0(x),  u_{t}(x,0)=\widetilde{u}_1(x),
\end{equation}
where
\begin{equation}\label{e5.4}
\widetilde{u}_0(x)=-x^2+x+1,
\widetilde{u}_1(x)=-\widetilde{u}_0(x), k(t)=\frac{1}{2}e^{-t},
\end{equation}
\begin{equation}\label{e5.5}
F(x,t)=(2-t)e^{-t}+U_{ex}(1-U_{ex}-U_{ex}^2),
\end{equation}
where
\begin{equation}\label{e5.6}
U_{ex}(x,t)=(-x^2+x+1)e^{-t}.
\end{equation}
The exact solution of the problem \eqref{e5.1}-\eqref{e5.3} with
$\widetilde{u}_0(x)$, $\widetilde{u}_1(x)$, $k(t)$ and $F(x,t)$
defined in \eqref{e5.4} and \eqref{e5.5} respectively, is the
function $U_{ex}$ given in \eqref{e5.6}. To solve problem
\eqref{e5.1}-\eqref{e5.3} numerically, we consider the
differential system for the unknowns $u_{j}(t)=u(x_{j},t)$,
$v_{j}(t)=\frac{du_j}{dt}(t)$, with $x_{j}=jh$, $h=\frac{1}{N}$,
$j=0,1,...,N $:
\begin{equation}\label{e5.7}
\begin{gathered}
\frac{du_{j}}{dt}(t)=v_{j}(t),  j=0,1,...,N,\\
\frac{dv_0}{dt}(t)=\frac{1}{h^2}\left[-(1+h)u_0(t)+u_1(t)\right]\\
-\frac{1}{h^2}\int_0^{t}k(t-s)\left[-(1+h)u_0(s)+u_1(s)\right]ds
-v_0^3(t)+u_0^2(t)+F(x_0,t),\\
\frac{dv_{j}}{dt}(t)=\frac{1}{h^2}\left[u_{j-1}(t)-2u_{j}(t)+u_{j+1}(t)\right]\\
-\frac{1}{h^2}\int_0^{t}k(t-s)\left[u_{j-1}(s)-2u_{j}(s)+u_{j+1}(s)\right]ds\\
-v_{j}^3(t)+u_{j}^2(t)+F(x_{j},t),  j=1,2,...,N-1,\\
\frac{dv_{N}}{dt}(t)=\frac{1}{h^2}\left[u_{N-1}(t)-(1+h)u_{N}(t)\right]\\
-\frac{1}{h^2}\int_0^{t}k(t-s)\left[u_{N-1}(s)-(1+h)u_{N}(s)\right]ds
-v_{N}^3(t)+u_{N}^2(t)+F(x_{N},t),\\
u_{j}(0)=\widetilde{u}_0(x_{j}), v_{j}(0)=\widetilde{u}_1(x_{j}),
j=0,1,...,N.
\end{gathered}
\end{equation}
To solve the nonlinear differential system \eqref{e5.7}, we use
the following linear recursive scheme generated by the nonlinear
term
\begin{equation}\label{e5.8}
\begin{gathered}
\frac{du_j^{(n)}}{dt}(t)=v_{j}^{(n)}(t), j=0,1,...,N,\\
\frac{dv_0^{(n)}}{dt}(t)=\frac{1}{h^2}\left[-(1+h)u_0^{(n)}(t)+u_1^{(n)}(t)\right]\\
-\frac{\Delta t}{h^2}\sum_{i=1}^{N_1-1}k(t-i\Delta
t)\left[-(1+h)u_0^{(n)}(i\Delta t)+u_1^{(n)}(i\Delta t)\right]\\
-\left(v_0^{(n-1)}(t)\right)^3+\left(u_0^{(n-1)}(t)\right)^3+F(x_0,t),\\
\frac{dv_j^{(n)}}{dt}(t)=\frac{1}{h^2}\left[u_{j-1}^{(n)}(t)-2u_{j}^{(n)}(t)+u_{j+1}^{(n)}(t)\right]\\
-\frac{\Delta t}{h^2}\sum_{i=1}^{N_1-1}k(t-i\Delta
t)\left[u_{j-1}^{(n)}(i\Delta t)-2u_{j}^{(n)}(i\Delta
t)+u_{j+1}^{(n)}(i\Delta t)\right]\\
-\left(v_{j}^{(n-1)}(t)\right)^3+\left(u_{j}^{(n-1)}(t)\right)^2+F(x_{j},t),
j=1,2,...,N-1,\\
\frac{dv_N^{(n)}}{dt}(t)=\frac{1}{h^2}\left[u_{N-1}^{(n)}(t)-(1+h)u_{N}^{(n)}(t)\right]\\
-\frac{\Delta t}{h^2}\sum_{i=1}^{N_1-1}k(t-i\Delta
t)\left[u_{N-1}^{(n)}(i\Delta t)-(1+h)u_{N}^{(n)}(i\Delta t)\right]\\
-\left(v_{N}^{(n-1)}(t)\right)^3+\left(u_{N}^{(n-1)}(t)\right)^2+F(x_{N},t),\\
u_{j}^{(n)}(0)=\widetilde{u}_0(x_{j}),
v_{j}^{(n)}(0)=\widetilde{u}_1(x_{j}), j=0,1,...,N,
\end{gathered}
\end{equation}
and where $u_{j}^{(n)}(i\Delta t)$, $i=1,...,N_1-1$,
$j=0,1,...,N$, of the system \eqref{e5.8} being calculated at the
time $t=N_1\Delta t$.

The latter system is solved by a spectral method and since the
matrix of this system is very ill-conditioned so we have to
regularize it by adding to the diagonal terms a small parameter in
order to have a good accuracy of the convergence.

\begin{center}
\includegraphics[width=7cm,height=7cm,
 keepaspectratio=true]{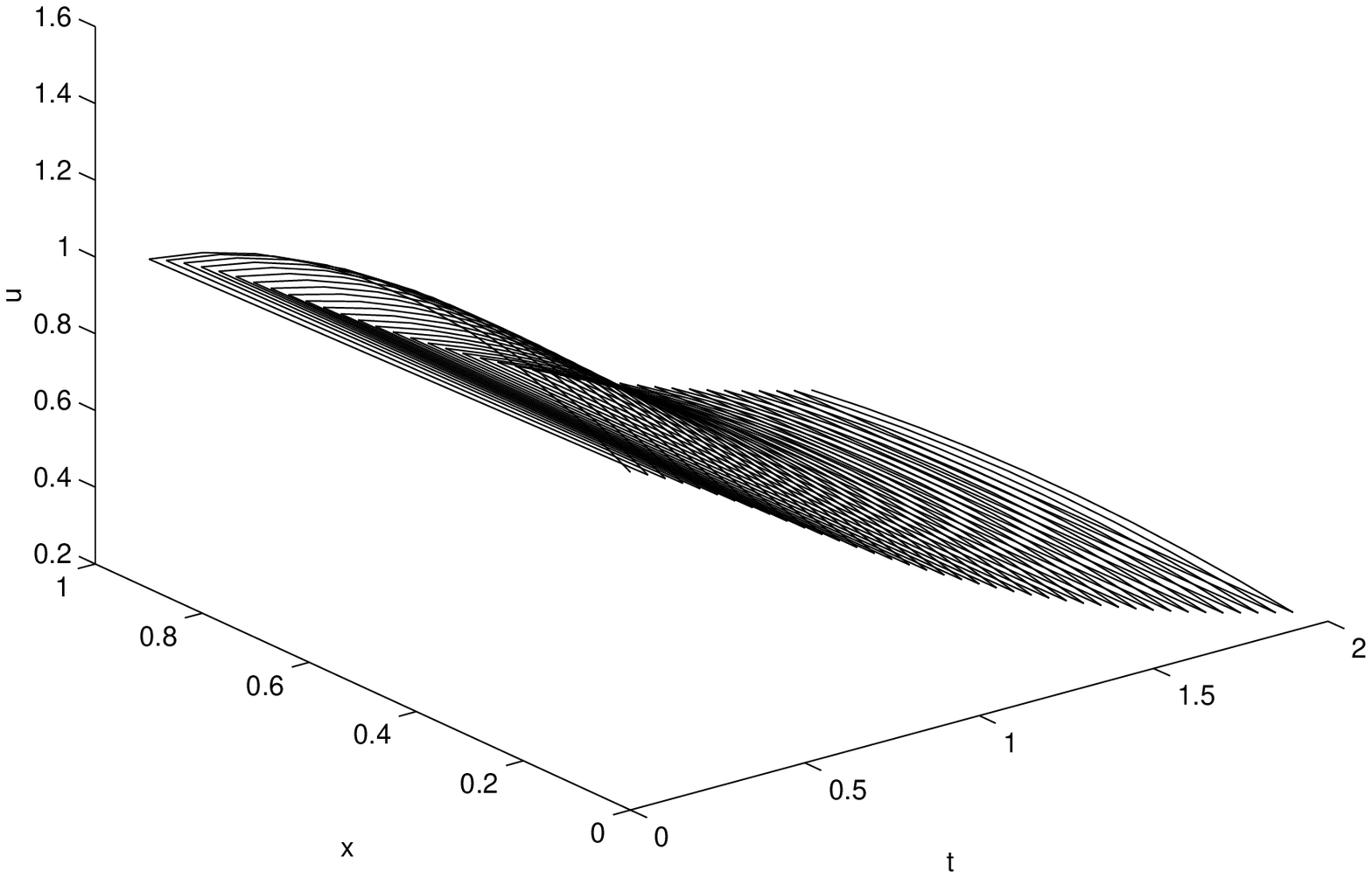}

Figure 1
\end{center}

In fig.1 we have drawn the approximated solution of the problem
\eqref{e5.1}-\eqref{e5.5} while fig.2 represents his corresponding
exact solution \eqref{e5.6}.

\begin{center}
\includegraphics[width=7cm,height=7cm,
keepaspectratio=true]{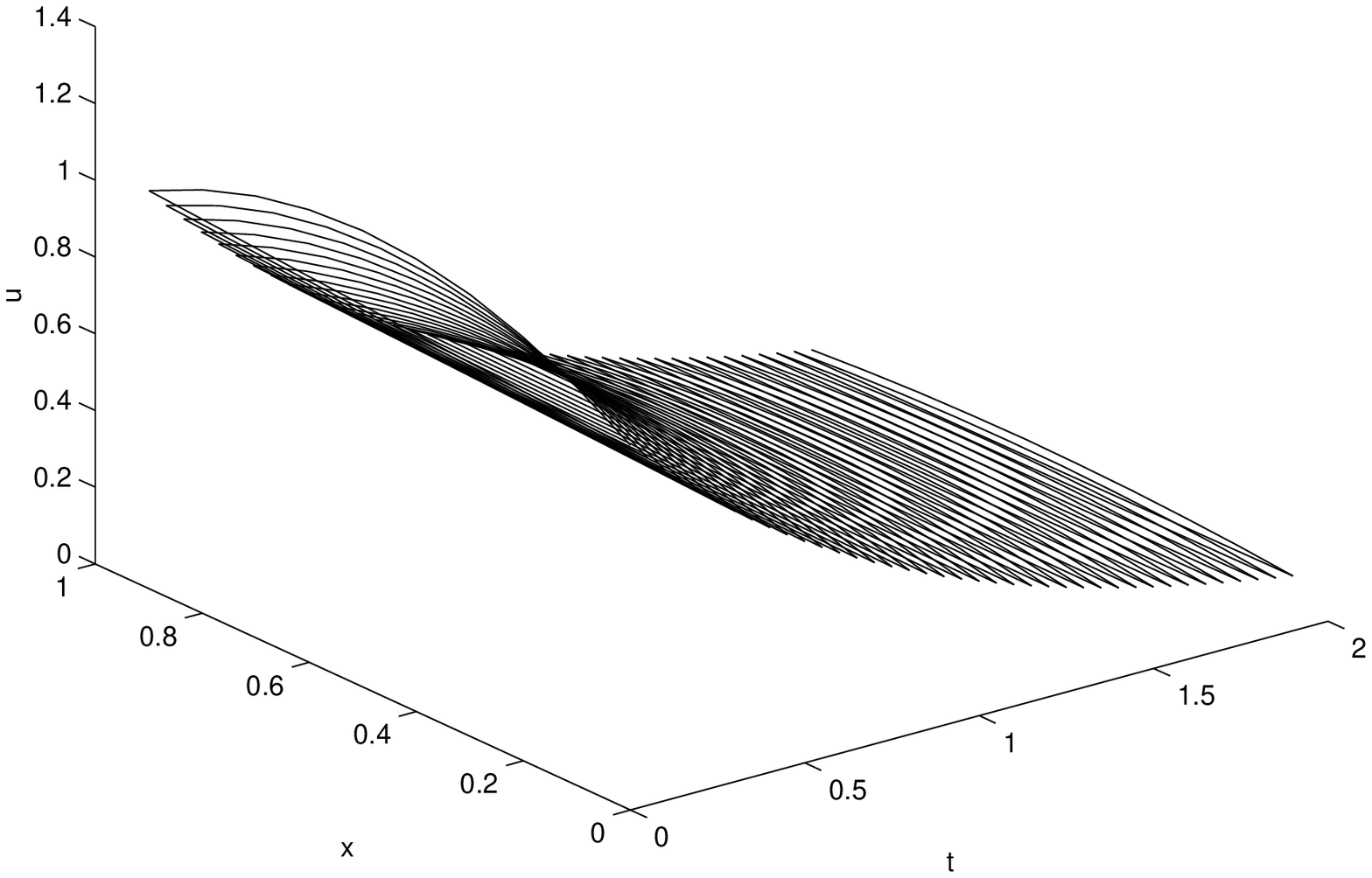}

Figure 2
\end{center}

The fig.3 corresponds to the surface $(x,t)\mapsto u(x,t)$
approximated solution in the case where $F(x,t)=0$. So in both cases
we notice the very good decay of these surfaces from  $T=0$ to
$T=2$. \\

\begin{center}
\includegraphics[width=7cm,height=7cm,
keepaspectratio=true]{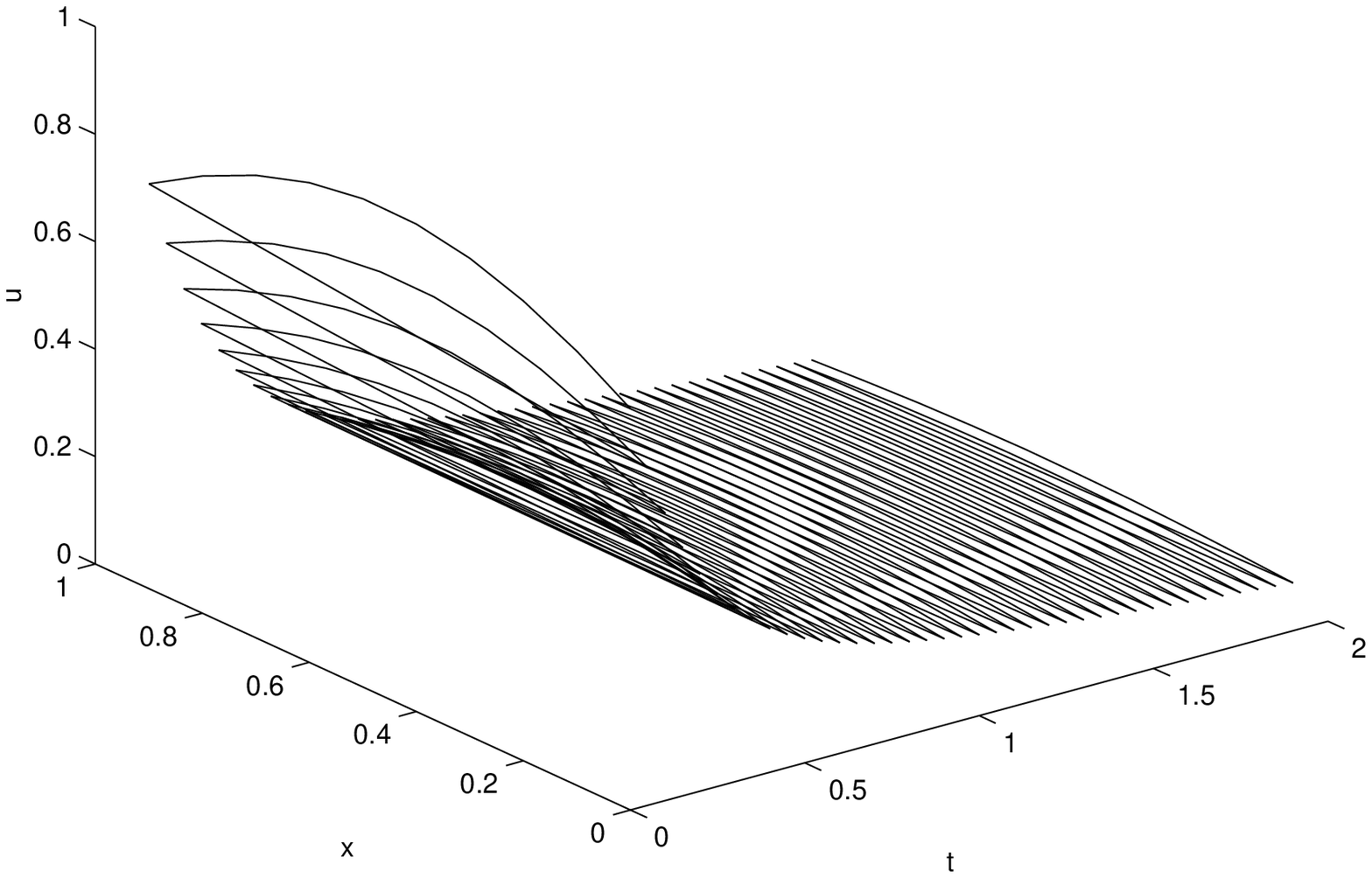}

Figure 3
\end{center}


\begin{thebibliography}{00}


\bibitem{b1} S. Berrimia, S. A. Messaoudi, \emph{Existence and decay of solutions of a viscoelastic
equation with a nonlinear source}, Nonlinear Anal. 64 (10) (2006)
2314-2331.

\bibitem{b2} M. Bergounioux, N. T. Long, A. P. N. Dinh, \emph{Mathematical model for a
shock problem involving a linear viscoelastic bar}, Nonlinear
Anal. {\bf 43} (2001), 547-561.

\bibitem{c1} M.M. Cavalcanti, V.N. Domingos Cavalcanti, J.A.
Soriano, \emph{Exponential decay for the solution of semilinear
viscoelastic wave equations with localized damping}, Electron. J.
Differential Equations 44 (2002) 1--14.

\bibitem{l1} J.L. Lions, W.A. Strauss, \emph{Some nonlinear evolution
equations}, Bull. Soc. Math., France, 93 (1965) 43-96.

\bibitem{l2} J. L. Lions, \emph{Quelques m\'{e}thodes de r\'{e}solution des probl\`{e}
mes aux limites nonlin\'{e}aires}, Dunod; Gauthier- Villars,
Paris. 1969.

\bibitem{l3} N.T. Long, A.P.N. Dinh, T.N. Diem, \emph{Linear recursive schemes and
asymptotic expansion associated with the Kirchhoff- Carrier
operator}, J. Math. Anal. Appl. 267 (2002) 116-134.

\bibitem{l4} N. T. Long, A. P. N. Dinh, T. N. Diem,
\emph{On a shock problem involving a nonlinear viscoelastic bar},
J. Boundary Value Problems, Hindawi Publishing Corporation, {\bf
2005} (3) (2005), 337-358.

\bibitem{l5} N. T. Long, L. V. Ut, N. T. T. Truc,
\emph{On a shock problem involving a linear viscoelastic bar},
Nonlinear Analysis, Theory, Methods \& Applications, Series A:
Theory and Methods, {\bf 63} (2) (2005), 198-224.

\bibitem{l6} N.T. Long, V.G. Giai, \emph{A wave equation
associated with mixed nonhomogeneous conditions: Global existence
and asymptotic expansion of solutions}, Nonlinear Analysis,
Theory, Methods and Applications, Series A: Theory and Methods, 66
(7) (2007), 1526-1546.

\bibitem{l7} N. T. Long, L. X. Truong,
\emph{Existence and asymptotic expansion for a viscoelastic
problem with a mixed nonhomogeneous condition}, Nonlinear
Analysis, Theory, Methods \& Applications, Series A: Theory and
Methods, (accepted for publication).

\end{thebibliography}
\end{document}